\documentclass[12pt]{article}
\usepackage{amssymb, amsmath, amsfonts}
\usepackage{mathrsfs}
\usepackage{dsfont}
\newcommand{\p}{\mathsf{P}}
\newcommand{\e}{\mathsf{E}}

\newcommand{\ent}{\widehat{H}_{n,k}}

\newcommand{\sgn}{\text{sgn}}

\newcommand{\ind}{{\mathbb{I}}}

\newtheorem{thm}{Theorem}

\newtheorem{cor}{Corollary}

\newtheorem{Lem}{Lemma}

\begin{document}

\begin{center}
{\Large{\bf Statistical Estimation of  Conditional Shannon Entropy }}
\end{center}

\begin{center}
{\large {\bf Alexander Bulinski\footnote{E-mail: bulinski@mech.math.msu.su}, Alexey Kozhevin\footnote{E-mail: kozhevin.alexey@gmail.com}}}
\end{center}
\begin{center}
{\it 
Dept. of Mathematics and Mechanics, Lomonosov Moscow State University\\ 
Steklov Mathematical Institute of Russian Academy of Sciences
%,\\ Moscow 119234, Russia
}
\end{center}
\vskip0.9cm
{\small
{\bf Abstract}
The new estimates of the conditional Shannon entropy are introduced in the framework
of the model describing a discrete response variable depending on a vector of $d$ factors
having a density w.r.t. the Lebesgue measure in $\mathbb{R}^d$.
Namely, the mixed-pair model $(X,Y)$ is considered where $X$ and $Y$ take values in $\mathbb{R}^d$ and an arbitrary finite set, respectively.
Such models include, for instance, the famous logistic regression. In contrast to the well-known Kozachenko -- Leonenko estimates of unconditional entropy the proposed estimates are
constructed by means of the certain spacial order statistics (or $k$-nearest neighbor statistics where $k=k_n$
depends on amount of observations $n$)   and a random number of i.i.d. observations
contained in the balls of specified random radii. The asymptotic unbiasedness and $L^2$-consistency of the new estimates are established under simple conditions. The obtained results can be
applied to the feature selection problem which is important, e.g., for medical and biological investigations.


\vskip1.7cm
\section{Introduction}
\vskip0.3cm
The entropy concept plays a prominent role in physics and mathematics,
see, e.g., \cite{Benguigui}.
On various approaches to the entropy definition we refer to the deep works by
L.Boltzmann, J.Gibbs, M.Plank, C.Shannon, A.N.Kolmogorov, Ya.G.Sinai, A.Renyi, C.Tsallis, A.S.Holevo. There are impor\-tant problems where one employs, for the due entropy, the statistical estimates
construct\-ed by means of i.i.d. observations.
For example such estimates are useful
in the feature selection theory (\cite{Peng}) and in the detection of  texture inhomogeneities (\cite{Alonso-Ruiz}).
We leave apart many other domains where the  entropies estimates are applied, see, e.g.,  \cite{Pal}.
There are a number of various approaches to the entropy estimation, we refer, e.g., to \cite{Archer}, \cite{Charzynska}, \cite{Hall},  \cite{Miller},
\cite{Paninski}, \cite{Sricharan}, \cite{Stowell}, \cite{Tsybakov}.

The main goal of the paper is to introduce the new statistical estimates of conditional Shannon entropy
for models where a discrete response variable, taking values in an arbitrary finite set, depends on a vector of factors (features) having
density w.r.t. the Lebesgue measure in $\mathbb{R}^d$. These models include the famous logistic regression (see, e.g., \cite{Hilbe}, \cite{Kleinbaum}).  The proposed estimates involve the  $k$-nearest neighbor statistics where $k=k_n$
depends on a number of observations $n$ (on the $k$-nearest neighbor statistics see, e.g.,
the recent book  \cite{Biau}). Note that our estimates do not employ the well-known Kozachenko-Leonenko statistics \cite{Kozachenko} used for estimation of the unconditional Shannon entropy. Under simple assumptions (cf., e.g., \cite{Berrett}, \cite{BulDim}, \cite{Delattre}, \cite{Gao}, \cite{Singh_1}) we establish the asymptotic unbiasedness
and $L^2$-consistency of our estimates when the sample size tends to infinity.
An interest in the study of conditional entropy is explained as follows. The mutual information
of two random vectors is represented by means of conditional entropy of one of them and unconditional entropy of another. That information characteristic of two random vectors facilitates the identification of relevant factors having impact on a response variable under consideration
(see, e.g., \cite{Bennasar},  \cite{Evans}, \cite{Fleuret}, \cite{Novovicova}, \cite{Vergara}).
Such analysis is useful in medical and biological studies. Thus statistical estimates of the mutual information involving new estimates will be
valuable for feature selection.

We stipulate that all the random variables and random vectors are defined on a probability space
$(\Omega,\mathcal{F},\p)$. Recall that the Shannon entropy (see \cite{Shannon}) of a discrete
random variable $Y$ taking values in a finite set $M$ with probabilities $P(y) := \p(Y = y)$, $y\in M$, and a (differential) entropy of a random vector $X$ in $\mathbb{R}^d$ having density
$f(\cdot)$, $x\in \mathbb{R}^d$, w.r.t. the Lebesgue measure $\mu$ are introduced by the following  respective formulas
\begin{equation}\label{Hdis}
H(Y) := -\e \log P(Y) = - \sum_{y \in M}  P(y)\log P(y),
\end{equation}
%\vspace{-0.2cm}
\begin{equation}\label{Hdif}
H(X) := -\e \log f(X) = - \int_{\mathbb{R}^d} f(x)\log f(x)  \, \mu(dx).
\end{equation}
Clearly, one can view the entropy as a function of a probability distribution since the above
formulas involve the laws of $X$ and $Y$.
Note that the probability distribution discretization  techniques for a random variable having a density (w.r.t. the Lebesgue measure) and evaluation of the Shannon entropy for thus arising random variables do not lead to the differential entropy as the mash of the discretization tends to zero
(see, e.g., Theorem 8.3.1 in \cite{Cover} and \cite{Montalvao}).
More generally, when a measure $\sigma$ is fixed on a measure space $(S,\mathcal{B})$, one can define the notion of the entropy of a probability measure $\nu$ given on the same space and
absolutely continuous w.r.t. $\sigma$. Namely,
whenever the following integral is well defined (and can take infinite values),
\begin{equation}\label{Hsigma}
H_{\sigma}(\nu):= - \int_S \log\left(\frac{d\nu}{d\sigma}\right)\,d\nu
\end{equation}
where $\frac{d\nu}{d\sigma}$ is the Radon -- Nikodym derivative.

If $Y$ has a law $\nu$ on $(M,2^M)$
then \eqref{Hdis} is a particular case of \eqref{Hsigma} where
$S=M$, $\mathcal{B}=2^M$ and  $\sigma$ is a counting measure on $M$.
If  $X$ has a law $\nu$ on  $(S,\mathcal{B})$ then \eqref{Hsigma}
leads to \eqref{Hdif} when $S=\mathbb{R}^d$, $\mathcal{B}=\mathcal{B}(\mathbb{R}^d)$ and $\sigma=\mu$.
The definition of the Kulback -- Leibler (see, e.g., \cite{Cover}, p.19, 251) relative entropy (or divergence) for two probability measures is closely related  to \eqref{Hsigma}.
We refer to \cite{Sason} where various kinds of $f$-divergences are compared.

Consider a random vector
$(X,Y)$ such that $X:\Omega\to \mathbb{R}^d$ ($d \in \mathbb{N}$) and $Y:\Omega \to M$.
Here $M$ is an arbitrary finite set. We assume that $\p(Y = y) > 0$ for each $y \in M$.
Suppose that there exists a measurable function
$f_{X,Y}: \mathbb{R}^d \times M \to \mathbb{R_{+}}$ such that, for any $B \in \mathcal{B}(\mathbb{R}^d)$ and  $y \in M$,
\begin{equation}\label{M}
\p(X \in B, Y = y) = \int_{B} f_{X,Y}(x,y) \, \mu(dx).
\end{equation}
In other words, $f_{X, Y}$ is a density of a random vector $(X,Y)$ w.r.t. measure $\sigma:=\mu \otimes \lambda$ on $\mathcal{B}(\mathbb{R}^d)\otimes 2^M$.
For $x\in \mathbb{R}^d$ and $y\in M$, let us define the following functions:
	\begin{gather*}
	f_{X}(x) := \sum_{y\in M}\int_{\mathbb{R}^d} f_{X,Y}(x,y) \, \mu(dx),\\
	f_{X|Y}(x|y) := \frac{f_{X,Y}(x,y)}{\p(Y=y)},
\end{gather*}
\begin{equation}\label{cd}
	f_{Y|X}(y|x) := \begin{cases}
			\frac{f_{X,Y}(x,y)}{f_X(x)}, & f_X(x) > 0, \\
			 0, & f_X(x) = 0. \\
		\end{cases}
\end{equation}	
Note that $f_X$ is a density of $X$, $f_{X|Y}$ is a conditional density of $X$ given $Y$, and $f_{Y|X}$ provides a conditional distribution of $Y$ given  $X$.
To simplify notation we will write $dx$ instead of $\mu(dx)$ and set $f(x,y) := f_{X,Y}(x,y)$, $f(y|x) := f_{Y|X}(y|x)$.

According to \eqref{Hsigma} (see also \cite{Nair}) the entropy of a vector $(X,Y)$ in the framework of model \eqref{M} is given by the formula
\begin{equation*}
H(X,Y) := -\e \log f(X,Y) = - \sum_{y\in M} \int_{\mathbb{R}^d} f(x,y) \log f(x,y) \, dx.
\end{equation*}
Introduce the conditional entropy of $Y$ given $X$
\begin{equation}\label{CH_1}
	H(Y|X) := -\e \log f(Y|X) = - \sum_{y\in M} \int_{\mathbb{R}^d} f(x,y) \log f(y|x) \, dx.
\end{equation}
One can verify that this conditional entropy $H(Y|X)$
is always finite.
Note that $H(Y|X) = H(X,Y) - H(X)$. The mutual information of $X$ and $Y$ is defined as
\begin{equation}\label{CH}
I(X,Y):=H(X)+H(Y)-H(X,Y)=H(X)-H(X|Y)=H(Y)-H(Y|X).
\end{equation}
It is well-known that $I(X,Y)\geq 0$. Moreover, $I(X,Y) = 0$ if
and only if $X$ and $Y$ are independent.
The latter statement is applied to the information approach for the identification
of relevant factors having an impact on a random response.
Mention in passing that extension of \eqref{CH} to the case of $n$ random vectors is fruitful as well (see, e.g., \cite{Doquire}).
There are a number of papers devoted to various estimates of (unconditional) entropy.
In this regard we indicate the recent work \cite{BulDim} where the estimates of
the Shannon differential entropy are studied and where one can find further references.

The scheme \eqref{M} under consideration comprises the famous logistic model widely used in the classification problems (see, e.g., \cite{Massaron}). Namely, let
$M=\{1,2\}$ and
\begin{equation}\label{L}
	\p(Y = 1|X = x) = \frac{1}{1+\exp\{-(w,x) - b\}}, \;\; x\in \mathbb{R}^d,\;\;w \in \mathbb{R}^d,\;\;b \in \mathbb{R},
\end{equation}
where $(\cdot,\cdot)$ is a scalar product in $\mathbb{R}^d$ and
	$\p(Y = 2|X = x) = 1 - \p(Y = 1|X = x)$.
Let $f_X$ be a vector $X$ density. Then
	\begin{gather*}
	f_{X,Y}(x,1) = \p(Y = 1|X = x) f_X(x), \\
	f_{X,Y}(x,2) = f_X(x) - f_{X,Y}(x,1).
	\end{gather*}
Note that there exist generalizations of logistic regression where a response variable
$Y$ takes more than two different values.

\vskip0.2cm
To conclude the introduction we mention that
in Section 2  statistical estimates of  $H(Y|X)$ are introduced and two principle results are formulated.
The proposed estimates are
constructed by means of the certain  $k$-nearest neighbor statistics (where $k=k_n$
depends on a number $n$ of i.i.d. observations $(X_1,Y_1),\ldots,(X_n,Y_n)$)   and a random number of observations
contained in the balls of specified random radii.
Under wide conditions the asymptotic unbiasedness and  $L^2$-consistency of our estimates are proved in Sections 4 and 5, respectively, whereas in Section 3 some auxiliary results are provided. Their proofs and
that of Corollary are given in Appendix.
The applications to the feature selection problems along with simulations will be considered separately. In particular, for considered vectors $(X,Y)$ our estimate of the conditional entropy
$H(Y|X)$  has advantages over estimates constructed as differences of statistical
estimates of $H(X,Y)$ and $H(X)$. Note also that other estimates of mutual information for discrete-continuous mixtures models
based on the Kraskov - St\"ogbauer - Grassberger \cite{Kraskov} approach were studied
in \cite{Coelho} and \cite{Gao} under different conditions. Also it is worth to emphasize that our estimates construction does not suppose the existence of any topological structure on
a set $M$ (thus we do not use the distances between $Y_i$ and $Y_j$, $i,j=1,\ldots,n$).

\vskip0.5cm

{\large {\bf 2. Main results}}
\vskip0.3cm

	Let $Z_1,Z_2,\ldots$ be a sequence of i.i.d. random vectors $Z_i = (X_i,Y_i)$, $i\in \mathbb{N}$, such that a distribution of
$Z_1$ coincides with one of the vector $(X,Y)$ described by model \eqref{M}. Introduce the estimate $H(Y|X)$  constructed by a sample $Z_1,\ldots,Z_n$ as follows
\begin{equation}\label{H}
	\widehat{H}_{n,k} = \frac{1}{n} \sum_{i=1}^n  \widehat{H}_{n,k,i}.
\end{equation}
Here $n\in\mathbb{N}$, $n>1$, $k=k(n) \in \{1,\dots,n-1\}$,
\begin{equation}\label{Hi}
	\widehat{H}_{n,k,i} = -\log (\xi_{n,k,i}(Z_1,\dots,Z_n)+1) + \log k,
\end{equation}
\begin{equation}\label{xi}
	\xi_{n,k,i}(Z_1,\dots,Z_n) := \sharp\{ j \in \{1,\dots,n\}\setminus \{i\} \colon Y_j = Y_i, \|X_i - X_j\| \leq \rho_{n,k,i}(X_1, \dots, X_n) \},
\end{equation}
$\sharp$ stands for the cardinality of a finite set, $\|\cdot\|$ is the Euclidean norm in $\mathbb{R}^d$ and
\begin{equation}\label{rodef}
\rho_{n,k,i}(X_1, \dots, X_n):= \|X_i - X_{i,(k)}\|,
\end{equation}
$X_{i,(k)}$ being the $k$-th nearest neighbor of $X_i$ in the sample $\{X_1, \dots, X_n\} \setminus \{X_i\}$
i.e. $\rho_{n,k,i}(X_1, \dots, X_n)$ is the Eulidean distance from $X_i$ to its $k$-th nearest neighbor. Clearly, the random variable $\xi_{n,k,i}(Z_1,\ldots,Z_n)$ takes values $0,1,\ldots,k$.
Observe that with probability one the points
$X_1,\ldots,X_n$ do not pair-wise coincide as the vector $X$
has a density.

Thus, in contrast to the well-known Kozachenko - Leonenko (\cite{Kozachenko}) estimate
of the Shannon differential entropy of a random vector, along with the distance to the
$k$-th nearest neighbor of $X_i$ in the sample $X_1,\ldots,X_n$ (without point $X_i$)
the principle role is played by random variables $\xi_{n,k,i}$, $i=1,\ldots,n$.
Namely, at first we find a random set $J\subset \{1,\ldots,n\}$, consisting of all the indexes
$j\in \{1,\ldots,n\}\setminus \{i\}$ such that
$X_j$ belongs to the ball
$B(X_i,\rho_{n,k,i})$ with a random center and a random radius. Then from the collection
$\{(X_j,Y_j),j\in J\}$ we take $\{(X_j,Y_j),j\in I_i\}$ where $I_i:= \{j\in J: Y_j=Y_i\}$. The collection of random variables $\{(X_j,Y_j),j\in I_i\}$ arises where $I_i$ is also a random set. The cardinality of this set
$I_i$,i.e. $\sharp I_i$, equals to the random variable  $\xi_{n,k,i}$, $i=1,\ldots,n$.

\vskip0.2cm
\noindent	
	\textbf{Definition}. A function $g \colon \mathbb{R}^d \to \mathbb{R}$ is called
locally constricted at a point $x$ in $\mathbb{R}^d$ if there exist strictly positive $R_0(x)$ and $C_0(x)$
such that
\begin{equation}\label{cw}
	\left|g(x) - \frac{1}{|B(x,R)|} \int_{B(x,R)} g(v)\, dv \right| \leq C_0(x)R\;\;\mbox{for}\;\;R\in (0,R_0(x))
\end{equation}
where $|B(x,R)|$ is a ball $B(x,R):=\{v\in \mathbb{R}^d:\|v-x\|\leq R\}$ volume, i.e. $|B(x,R)|=\mu(B(x,R))$.
	A function $g$ is $C_0$-constricted if it is locally constricted
for $\mu$-almost all points $x \in \mathbb{R}^d$ and, moreover, for such $x$ one has $C_0(x)\leq C_0$ and  $R_0(x)\geq R_0$ where $C_0$ and $R_0$ are strictly positive constants.
\vskip0.2cm	
\noindent	
	\textbf{Remark 1}. If a function $g:\mathbb{R}^d\to \mathbb{R}$
satisfies the Lipschitz condition at $x\in \mathbb{R}^d$ with a factor $C(x)$, that is $|g(v)-g(x)|\leq C(x)\|v-x\|$ for all $v\in \mathbb{R}^d$, then  \eqref{cw} is valid for any $R_0(x)>0$.
It is easily seen that if $g(x)$, $x\in \mathbb{R}^d$, is a density of non-degenerate Gaussian law
then this function is $C_0$-constricted.

	\begin{thm}\label{th_1}
Let in the framework of model \eqref{M} the following conditions be satisfied.
For each fixed $y\in M$ and $\mu$-almost all $x\in \mathbb{R}^d$, a function $f(x,y)$, i.e. $f(\cdot,y)$, is strictly
positive and $C_0$-constricted,
\begin{equation}\label{k_cond}
			k = k_n \propto n^{\alpha}
\end{equation}
for some $\alpha \in (0, 1)$, and, for some $\varepsilon > 0$,
\begin{equation}\label{e_log}
			\e |\log f_X(X)|^{1+\varepsilon}  < \infty
\end{equation}
where $f_X(\cdot)$ is a density of $X$.
Then
\begin{equation}\label{aunb}
		\e \ent \to H(Y|X),\;\;n\to \infty,
\end{equation}		
i.e. $\ent$ is an asymptotically unbiased estimate of $H(Y|X)$.
	\end{thm}

	\begin{thm}\label{th_2}
Let the condition \eqref{e_log} of Theorem \ref{th_1} be replaced by the following one. For some  $\varepsilon > 0$,
\begin{equation}\label{AB}
		\e |\log f_X(X)|^{2+\varepsilon}  < \infty.
\end{equation}
Then
		$$
		\e(\ent - H(Y|X))^2 \to 0\;\;\mbox{as}\;\;n \to \infty,
		$$
i.e. $\ent$ is an $L_2$-consistent estimate of $H(Y|X)$.
	\end{thm}

\begin{cor}\label{Gauss}
Let in the framework of model \eqref{M}, for each $y\in M$, the function $f(\cdot,y)$
 be a density of non-degenerate Gaussian law in $\mathbb{R}^d$ $($with mean vector and covariance matrix depending on $y)$. Then $\widehat{H}_{n,k}$, where $k(n)$ satisfies \eqref{k_cond}, are asymptotically unbiased and  $L^2$-consistent estimates of $H(Y|X)$ as
$n\to \infty$.

\end{cor}

\vskip0.5cm
{\large {\bf 3. Auxiliary results}}
\vskip0.3cm
In this Section, as previously, we consider i.i.d. vectors
$Z_1,Z_2,\ldots$, having the same distribution as $(X,Y)$ in the framework of model \eqref{M}.
Using notation introduced in Section 2,  for $x\in \mathbb{R}^d$, $y\in M$, $k\in\{1,\ldots,n-1\}$ and $n>2$, set
$$
\rho_{n,k,1}(x):= \rho_{n,k,1}(x,X_2,\ldots,X_n),\;\;
\xi_{n,k,1}(x,y):= \xi_{n,k,1}((x,y),Z_2,\ldots,Z_n).
$$
These random variables depend, respectively, on $x,X_2,\ldots,X_n$ and $(x,y),Z_2,\ldots,Z_n$.
To simplify notation we omit the random arguments of these functions.

Now we formulate two auxiliary results concerning conditional distributions of random variables
playing an essential role in the asymptotical behavior analysis of the conditional entropy estimates.
It turns out surprisingly that the mentioned conditional distributions
are specified mixtures of certain binomial laws with explicitly indicated weight coefficients.
The proofs of these results are provided in Appendix. We write $\p_{\eta}$ for distribution of
a random vector (or variable) $\eta$.
\begin{Lem}\label{lem_1}
For any $y\in M$, $x\in \mathbb{R}^d$,  $r=0,1,\ldots,k$ and
$\p_{\rho_{n,k,1}(x)}$-almost all $t\in (0,\infty)$, the following relation holds:
$$
\p(\xi_{n,k,1}(x,y)=r|\rho_{n,k,1}(x)=t)
$$
$$
= \binom{k-1}{r}  \p(Y=y|\|X-x\|\leq t)^r(1-\p(Y=y|\|X-x\|\leq t)^{k-1-r}\alpha(x,y,t)
$$
$$
+ \binom{k-1}{r-1}\p(Y=y|\|X-x\|\leq t)^{r-1}(1-\p(Y=y|\|X-x\|\leq t)^{k-r}(1-\alpha(x,y,t))
$$
where $
\alpha(x,y,t)= \p(Y\neq y|\|X-x\|=t)
$ and
$\binom{N}{m}:=0$ for $m<0$ and $m>N$ $(N\in \mathbb{N}$, $m\in \mathbb{Z})$.
\end{Lem}
\noindent
{\bf Remark 2}.
As usual, for random vectors  $\eta:\Omega\to \mathbb{R}^q$,
$\zeta:\Omega\to \mathbb{R}^s$, and for $B\in \mathcal{B}(\mathbb{R}^q)$, $x\in \mathbb{R}^s$, the notation $\p(\eta \in B|\zeta=x)=\varphi(x)$ means that one takes a Borel function $\varphi(x)$, $x\in \mathbb{R}^s$, such that $\p(\eta \in B|\zeta)=\varphi(\zeta)$. The function $\varphi$ is defined uniquely  $\p_{\zeta}$-almost sure, see, e.g., \cite{Shiryaev}, v.1, Ch.II, Section 5.
\vskip0.2cm

Let $z_j=(x_j,y_j)$ where $x_j\in \mathbb{R}^d$, $y_j\in M$, $j\in \{1,2\}$.
For $n>2$, let us define the random variables
	\begin{equation}\label{ro}
\rho_{n,k,j}(x_1, x_2) := \rho_{n,k,j}(x_1 ,x_2, X_3, \ldots,X_n),
\end{equation}
\begin{equation}\label{hi}
	\xi_{n,k,j}(z_1,z_2) := \xi_{n,k,j}(z_1, z_2, Z_3, \dots, Z_n).
	\end{equation}
Again we omit the random arguments of these functions.
	\begin{Lem}\label{lem_2}
Let $z_j=(x_j,y_j)$ where $x_j\in \mathbb{R}^d$, $y_j\in M$, $j=1,2$, $x_1\neq x_2$.
Introduce a random vector $\zeta:=(\rho_{n,k,1}(x_1, x_2),\rho_{n,k,2}(x_1, x_2))$.
Then, for any $n>2$, $k\in \{0,1,\ldots,[(n-2)/2]\}$, where $[\cdot]$ is the integer part of a number,  any $r_1,r_2\in \{0, \dots, k\}$ and
$\p_{\zeta}$-almost all $(t_1,t_2)$ such that $t_1>0$, $t_2>0$ and $t_1, t_2 < |x_1 - x_2|/2$ the following relation is valid:
\begin{equation}\label{conind}
\p(\xi_{n,k,1}(z_1,z_2)=r_1,\xi_{n,k,2}(z_1,z_2)=r_2|\zeta=(t_1,t_2))
=\prod_{j=1}^2 \p(\xi_{n,k,j}(z_1,z_2)=r_j|\zeta=(t_1,t_2)).
\end{equation}
Moreover,
$$
\p\left(\xi_{n,k,j}(z_1,z_2) = r_j|\zeta=(t_1,t_2)\right)
$$
\begin{equation}\label{claw}		
= \binom{k-1}{r_j}
		p_j^{r_j}
		(1 - p_j)^{k-1-r_j}\alpha(x_j,y_j,t_j)
%\p(Y\neq y_j|\|X-x_j\|=t_j)
+ \binom{k-1}{r_j-1}
		p_j^{r_j-1}
		(1-p_j)^{k-r_j}
%\p(Y = y_j | \|X-x_j\|=t_j)
(1-\alpha(x_j,y_j,t_j))
		\end{equation}		
where $p_j = \p(Y = y_j|X \in B(x_j, t_j))$,  $j=1,2$, and $\alpha(x,y,t)$ is the same as
in Lemma \ref{lem_1}.
\end{Lem}
We will also employ the following elementary results.
\begin{Lem}\label{lem_3}
Let $W$ be a random variable having finite  $\e W$, and $V$ be a random vector with values
in $\mathbb{R}^m$ such that $\p(V \in B) > 0$ where $B \in \mathcal{B}(\mathbb{R}^m)$.
Then
	$$
	\e (W|V \in B) = \int\limits_{B} \e(W|V=y) \, \widetilde{\p}_{V,B}(dy),
	$$
where $\e (W|V \in B):=\frac{1}{\p(V\in B)}\e (W\ind\{V\in B\})$
and $\widetilde{\p}_{V,B}(A):= \p (V\in A|V\in B)$, $A \in \mathcal{B}(\mathbb{R}^m)$.
\end{Lem}

\begin{Lem}\label{lem_4}
Let $\xi$, $\eta$ be some random variables and $\e |\xi| < \infty$. Assume that a random variable $\zeta$
takes values in a finite or countable set $S$. Then, for $\p_\eta$-almost all $t$, one has
	$$
	\e(\xi|\eta = t) = \sum_{r \in S} \e (\xi|\zeta=r, \eta=t) \p(\zeta=r|\eta=t).
	$$
\end{Lem}

\vskip0.3cm
{\large {\bf 4. Proof of Theorem \ref{th_1}}}	
\vskip0.2cm	
\noindent
Observations $Z_1,Z_2,\ldots$ have identical distribution. Thus
$
\e \widehat{H}_{n,k} = -\e \log \left(\frac{\xi_{n,k,1}+1}{k}\right)
$
and one has to prove that
$$
-\e \log \left(\frac{\xi_{n,k,1}+1}{k}\right)\to H(Y|X),\;\;n\to \infty.
$$	
Taking into account the independence of $Z_1$ and $Z_2,\ldots,Z_n$,  we get
$$
\e\left(\log\left(\frac{\xi_{n,k,1}+1}{k}\right)\bigg|X_1=x,Y_1=y\right)
=\e\log\left(\frac{\xi_{n,k,1}(x,y)+1}{k}\right)
:=h_{n,k}(x,y).
$$
Consequently,
$$
\e \log\!\left(\!\frac{\xi_{n,k,1}+1}{k}\right)
= \sum_{y\in M}\int_{\mathbb{R}^d}h_{n,k}(x,y)f(x,y)dx.
$$

Introduce parameters $\theta,
\nu > 0$. In the sequel we will make an appropriate choice of these parameters. Due to \eqref{e_log}, for $n>2$ and $y\in M$, we come to relations
\begin{gather}
\int\limits_{\{x: f_X(x) \leq n^{-\theta}\}}\!\!\!\! f_X(x) \, dx
\leq \int\limits_{\{x: f_X(x) \leq n^{-\theta}\}}\!\! \frac{|\log f_X(x)|^{1+\varepsilon}}{|\log n^{-\theta}|^{1+\varepsilon}} f_X(x) \, dx
\leq \frac{1}{(\theta \log n)^{1+\varepsilon}} \e |\log f_X(X)|^{1+\varepsilon},\label{dens_1}\\
\int\limits_{\{x: f(y|x) \leq n^{-\nu}\}} f(x,y) \, dx =  \int\limits_{\{x: f(y|x) \leq n^{-\nu}\}} f(y|x) f_X(x) \, dx
\leq  n^{-\nu} \int\limits_{\mathbb{R}^d} f_X(x) \, dx
=  n^{-\nu}. \label{dens_3}
\end{gather}
For $n\in \mathbb{N}$, take
$\theta_n := n^{-\theta}$, $\nu_n := n^{-\nu}$ and consider the sets
\begin{equation}\label{B1_B2}
B_{1,n} := \bigcap_{y \in M} \{x \in \mathbb{R}^d \colon f(y|x) > \nu_n\} \cap \{x \in \mathbb{R}^d \colon   f_X(x) > \theta_n\},\;\;
	B_{2,n} := \mathbb{R}^d \setminus B_{1,n}.
\end{equation}
One can write
$
\e \log\left(\frac{\xi_{n,k,1}+1}{k}\right) = I_1(n,k) + I_2(n,k)
$
where
\begin{gather*}
	I_j(n,k) := \sum_{y\in M} \int\limits_{B_{j,n}}h_{n,k}(x,y)f(x,y)dx,\;\;j=1,2.
\end{gather*}
 For $k>1$, all $x\in \mathbb{R}^d$ and $y\in M$, the inequality $|h_{n,k}(x,y)|\leq \log k$ is valid because
$\xi_{n,k,1}(x,y)$ takes values $0,1,\ldots,k$. Thus
	\begin{equation}\label{bound}
	|I_2(n,k)| \leq \log k \left(\,\int\limits_{\left\{x: f_X(x) \leq \theta_n \right\} } f_X(x) \, dx
+ \sum_{y \in M} \int\limits_{\left\{x: f(y|x) \leq \nu_n \right\} } f(x,y) \, dx \,\right).
	\end{equation}
According to \eqref{dens_1} and \eqref{dens_3}
we infer that $I_2(n,k) \to 0$, $n\to \infty$, since $k \propto n^\alpha$.

Fix parameter $\beta >0$ and note that
\begin{gather*}
\e \left(\log\left(\frac{\xi_{n,k,1}(x,y)+1}{k}\right)\Bigg|n^{\beta}\rho_{n,k,1}(x)\right)
= \sum_{r=0}^k \log\left(\frac{r+1}{k}\right)\p(\xi_{n,k,1}(x,y)=r|n^{\beta}\rho_{n,k,1}(x)).
\end{gather*}
Hence,
$$
h_{n,k}(x,y)= \sum_{r=0}^{k} \int_{(0,\infty)}\log\left(\frac{r+1}{k}\right)\p(\xi_{n,k,1}(x,y)=r|n^{\beta}\rho_{n,k,1}(x)=u)
f^{(k)}_{n,x,\beta}(u)du
$$
where $f^{(k)}_{n,x,\beta}(\cdot)$ is a density of a positive random variable $n^{\beta}\rho_{n,k,1}(x)$,
and a density of a random variable $\rho_{n,k,1}(x)$ is indicated in the proof of Lemma 1
(see Appendix, formula \eqref{dens}).

Now we fix an arbitrary  $\delta >0$ and write
$I_1(n,k)= S_1(n,k)+S_2(n,k)$ where, for  $V_1=(\delta,\infty)$,  $V_2=(0,\delta]$ and
$j=1,2$,
\begin{gather*}
	S_j(n,k) \!\!:=\!   \sum_{y\in M}\!\int_{B_{1,n}}\!\sum_{r=0}^k \!\int_{V_j}\!\!\!\!\log\left(\frac{r+1}{k}\right)\p(\xi_{n,k,1}(x,y)\!=\!r|\rho_{n,k,1}(x)\!=\!un^{-\beta})
f^{(k)}_{n,x,\beta}(u)du
f(x,y)dx.
\end{gather*}

The rest of the proof is divided into two steps.

Step 1. Let us show that
$S_1(n,k) \to 0$, $n\to \infty$.
We find an upper bound for $|S_1(n,k)|$.
For all $n\in \mathbb{N}$, $1\leq k\leq n$ and $x,y\in \mathbb{R}^d$, the variable $\xi_{n,k,1}(x,y)$ takes values  $0,\ldots,k$. Therefore, for $k>1$, one has
$$
\sum_{r=0}^k \left|\log\left(\frac{r+1}{k}\right)\right|\p(\xi_{n,k,1}(x,y)=r|\rho_{n,k,1}(x)=un^{-\beta})
\leq \log k.
$$
Thus
$$
	|S_1(n,k)| \leq \log k  \sum_{y=1}^{m} \int\limits_{B_{1,n}} \int_{(\delta,\infty)} f^{(k)}_{n,x,\beta}(u)f(x,y) \, du \, dx
$$
\vspace{-0.3cm}
\begin{gather}\label{ineq1}
	= \log k   \int\limits_{B_{1,n}} \p\left(\rho_{n,k,1}(x) > \delta n^{-\beta} \right) f_X(x) \, dx = \log k \,
\int\limits_{B_{1,n}} \p(\eta_n(\beta,\delta,x) \leq k-1) f_X(x)\, dx,
\end{gather}
here $	\eta_n(\beta,u,x) \sim Bin(n-1, p_n(\beta,u,x))$, i.e. $\eta_n(\beta,u,x)$ has a binomial law with parameters $n-1$ and $p_n(\beta,u,x)$ where
	\begin{gather*}
	{p}_n(\beta,u,x) := \p(X \in U_n(\beta,u,x)) = \int_{U_n(\beta,u,x)} f_X(v) \, dv, \\
	U_n(\beta,u,x) = \{v \in \mathbb{R}^d \colon \|v - x\| \leq  u n^{-\beta} \},\;\;u>0, \;\; x \in \mathbb{R}^d.
	\end{gather*}
Indeed, as $f^{(k)}_{n,x,\beta}(\cdot)$ is a density of a variable $\rho_{n,k}(x) n^{\beta}$, we
can write
	\begin{equation}\label{f1}
	\int_{(\delta,\infty)} f^{(k)}_{n,x,\beta}(u) \, du
= \p(\rho_{n,k,1}(x) > \delta n^{-\beta}).
	\end{equation}
The event
	$\{\omega: \rho_{n,k,1}(x) > \delta n^{-\beta}  \}$
means that in a  ball  $B(x,\delta n^{-\beta})$
one can find no more than $k-1$ point among $\{X_i\}_{i=2}^n$.
The independence of the observations yields
$	\p(\rho_{n,k,1}(x) > \delta n^{-\beta}) = \p(\eta_n(\beta,\delta,x) \leq k-1)
$.
According to the inequality for binomial sums proved in \cite{Zub},
for any $n>1$, $k=0,1,\ldots,n-1$ and all considered values  $\beta,\delta$ and $x$,   the
following bound holds
$$
\p(\eta_n(\beta,\delta,x) \leq k-1)
$$
\begin{equation}\label{ineq2}
\leq \Phi\left(\sgn\left(\frac{k}{n-1} - p_n(\beta,\delta,x)\right) \sqrt{2(n-1)  h \left( \frac{k}{n-1}, p_n(\beta,\delta,x)\right)} \right)
\end{equation}
where $\Phi(\cdot)$ is the
distribution function of a standard normal random variable,
	\begin{gather*}
	\sgn(t) = \begin{cases}
	\;\;1, & t > 0, \\
	\;\;0, & t = 0, \\
	\!-1, & t < 0,
	\end{cases} \\
	h(t,s) = t\log\left(\frac{t}{s}\right) + (1-t)\log\left(\frac{1-t}{1-s}\right),\; \;
s,t\in (0,1).
	\end{gather*}
For each $y\in M$, the function $f(\cdot,y)$ is $C_0$-constricted, therefore,
for $\mu$-almost all  $x\in \mathbb{R}^d$, each $u\in (0,\delta]$ and any
$n$ large enough,
\begin{equation}\label{d_1}
\left|\frac{\p(X\in U_n(\beta,u,x),Y=y)}{|U_n(\beta,u,x)|}- f(x,y)\right|
\leq C_0u n^{-\beta}.
\end{equation}
Since $f_X(x)=\sum_{y\in M} f(x,y)$, for
$x$, $u$ and $\beta$ under consideration, we get
\begin{equation}\label{e_imp}
\left|f_X(x)- \frac{p_n(\beta,u,x)}{|U_n(\beta,u,x)|}\right| \leq \sharp M C_0 u n^{-\beta}= C u n^{-\beta}
\end{equation}
where $C=\sharp M C_0$.
This implies that,
for arbitrary $\delta > 0$, $\beta >0$,
$\mu$-almost all $x\in B_{1,n}$ and for any $n\geq N_0$, where $N_0=N_0(\delta,\beta)$,
the following inequality is satisfied
\begin{equation}\label{est_lower}
p_n(\beta,\delta,x) \geq V_d \delta^d n^{-d\beta} (f_X(x) - C\delta n^{-\beta})) \geq V_d \delta^d n^{-d\beta}(\theta_n - C\delta n^{-\beta}).
\end{equation}
Recall that $V_d=|B(0,1)|$, $B(0,1)\subset \mathbb{R}^d$.
Now we can obtain the upper bound for the argument of a function $\sgn$ in formula \eqref{ineq2}. In view of \eqref{est_lower} one has
	\begin{gather*}
	\frac{k}{n-1} - p_n(\beta,\delta,x) \leq \frac{k}{n-1} - V_d \delta^d n^{-d\beta}(\theta_n - C\delta n^{-\beta}).
	\end{gather*}
Take $\beta >\theta$. Then
$n^{-\beta} = o(\theta_n)$ as $n\to \infty$.
According to \eqref{k_cond} we get $\frac{k}{n-1} \propto n^{\alpha-1}$, $n\to \infty$.
Let parameters $%\alpha,
\beta$ and $\theta$ be such that
\begin{equation}\label{c_1}
\alpha -1 <-d\beta -\theta.
\end{equation}
For \eqref{c_1} validity it is sufficient that $(d+1)\theta < 1-\alpha$, because
we can choose $\beta> \theta$ arbitrary close to $\theta$.
Then $n^{\alpha-1} = o(n^{-d\beta-\theta})$, $n\to \infty$. Thus there exists $N \in \mathbb{N}$ (where $N=N(C,d,\alpha,\delta,\beta,\theta)$) such that if  $n > N$, then $	\frac{k}{n-1} - p_n(\beta,\delta,x) < 0$ for $\mu$-almost all $x \in B_{1,n}$, which yields
\begin{gather}
	\sgn\left(\frac{k}{n-1} - p_n(\beta,\delta, x)\right) = -1. \label{sgn_ineq}
	\end{gather}
For $s,t \in (0,1)$ introduce the functions
$$
L_{1,n}(t,s)= 2(n-1)(1-t)\log\left(\frac{1-t}{1-s}\right),\;\;
L_{2,n}(t,s)= 2(n-1)t\log \left(\frac{t}{s}\right).
$$
Then $2(n-1)h(t,s)=L_{1,n}(t,s)+L_{2,n}(t,s)$.
Now, for $t=\frac{k}{n-1}$ and $s=p_n(\beta,\delta,x)$ consider the behavior of the functions
$L_{1,n}(t,s)$ and $L_{2,n}(t,s)$ (for $\mu$-almost all $x\in B_{1,n}$) as $n\to \infty$.
Applying \eqref{est_lower}, for $n\geq N_0$, we come to the bound
$$
L_{1,n}(t,s) = 2(n-1-k)\left(\log\left(1-\frac{k}{n-1}\right)-\log(1-p_n(\beta,\delta,x))\right)
$$
\begin{equation}\label{L2}
\geq 2(n-1-k)\left(\log\left(1-\frac{k}{n-1}\right)-
\log(1- V_d\delta^dn^{-d\beta}(\theta_n - C\delta n^{-\beta}))\right):= \mathcal{L}_{1,n}(k).
\end{equation}
Evidently, $\mathcal{L}_{1,n}(k)$ depends not only on $n$ and $k$, but also
on a collection of parameters appearing in â \eqref{L2}.
Note that $\log(1+z) = z + o(z)$ 	as $z \to 0$. Hence,
in view of \eqref{c_1} and since $\beta > \theta$ we get
$\mathcal{L}_{1,n}(k) \propto n^{1-d\beta-\theta}$, as $n\to \infty$.
For the same $t$, $s$ and $x\in B_{1,n}$, taking into account that $0< p_n(\beta,\delta,x)\leq 1$,
we obtain
$$
L_{2,n}(t,s)= 2k \left(\log \left(\frac{k}{n-1}\right)- \log p_n(\beta,\delta,x)\right)
\geq 2k \log \left(\frac{k}{n-1}\right):= \mathcal{L}_{2,n}(k).
$$
Therefore,
$
\mathcal{L}_{2,n}(k)\propto n^{\alpha}\log n
$ ïðè $n\to \infty$.
Thus according to
\eqref{c_1}
we conclude that, for all $n$ large enough, and for $\mu$-almost all $x\in B_{1,n}$,
\begin{equation}\label{R_cond}
	2(n-1)H\left( \frac{k}{n-1}, p_n(\beta,\delta,x)\right) \geq \mathcal{L}_{1,n}(k) + \mathcal{L}_{2,n}(k) := R_n \propto
	n^{1-d\beta-\theta}
\end{equation}
where $k=k_n\propto n^{\alpha}$. Here $R_n$ depends not only on $n$, but also on  $\alpha,\beta$, $\delta$, $\theta$ and $d$.
Therefore, \eqref{R_cond} gives, for all $n$ large enough, an estimate
$$
	\Phi\left(\sgn\left(\frac{k}{n-1} - p_n(\beta,\delta,x)\right) \sqrt{2(n-1)  h \left( \frac{k}{n-1}, p_n(\beta,\delta,x)\right)} \right) \leq \Phi(-\sqrt{R_n}).
$$
Since $\Phi(-z) \leq \frac{1}{\sqrt{2 \pi}z} e^{-\frac{z^2}{2}}$ for $z>1$,
taking into account  inequalities  \eqref{ineq1}, \eqref{ineq2} and \eqref{R_cond}
we get, for $k=k_n \propto n^{\alpha}$,
$x\in B_{1,n}$
and all  $n$ large enough,
\begin{equation}\label{ineq4}
\int_{(\delta,\infty)} f^{(k)}_{n,x,\beta}(u)\, du \leq \frac{1}{\sqrt{2 \pi R_n}} e^{-R_n/2},
\end{equation}
$$
|S_1(n,k)| \leq \frac{\log k_n}{\sqrt{2 \pi R_n}} e^{-R_n/2} \to 0,\;\;n\to \infty.
$$
Hence the proof of Step 1 is complete.

Step 2. We show that $S_2(n,k)\to - H(Y|X)$, $n\to \infty$.
For $\p_{n^{\beta}\rho_{n,k,1}(x)}$-almost all $u$ by virtue of Lemma \ref{lem_1}
$$
	\e \left( \log  (\xi_{n,k,1}(x,y)+1)  |\rho_{n,k,1}(x) = un^{-\beta} \right)
$$
\begin{equation}\label{st_2}
	= \p(Y \ne y | \|X - x\| =
un^{-\beta}) \e\log(\mu_{n}+1)   + \p(Y = y | \|X - x\| = un^{-\beta})\e\log(\mu_{n}+2)
\end{equation}
where the random variable $\mu_n$ does not depend on $(X,Y)$ and
\begin{equation}\label{mu_n}
	\mu_{n}=\mu_n(k,\beta,u,x,y) \sim Bin(k-1,P_n(\beta,u, x,y)),
\end{equation}
\begin{equation}\label{Pn}
	P_n(\beta,u, x,y) := \p(Y=y| X \in U_n(\beta,u, x)).
\end{equation}
Note that $\p(X\in B(x,t))>0$ for each $x\in \mathbb{R}^d$ and any $t>0$ since $f_X(\cdot)$ is strictly positive $\mu$-almost everywhere.
At first we study, for $u\in (0,\delta]$,  $y\in M$ and $\mu$-almost all
$x\in \mathbb{R}^d$ (such that $f_X(x)>0$), the convergence rate of $P_n(\beta,u, x,y)$ to $f(y|x)$ as $n\to \infty$. One has
$$
|P_n(\beta,u, x,y)-f(y|x)|
= \frac{|U_n(\beta,u,x)|}{\p(X\in U_n(\beta,u,x))}\Bigg|\Bigg\{
\left(
\frac{\p(X\in U_n(\beta,u,x),Y=y)}{|U_n(\beta,u,x)|}- f(x,y)\right)
$$
$$
+ f(y|x)
\left(f_X(x)- \frac{\p(X\in U_n(\beta,u,x))}{|U_n(\beta,u,x)|}\right)
\Bigg\}\Bigg|.
$$
For all $y\in M$ and $\mu$-almost all $x\in \mathbb{R}^d$,
the Lebesgue theorem on measures differentiation (see, e.g., Theorem 25.17 \cite{Yeh}) gives that
$$
f(y|x) =\lim_{R\to 0+}\frac{\p(Y=y,X\in B(x,R))}{\p(X\in B(x,R))}\leq 1.
$$
Hence, due to \eqref{d_1} and \eqref{e_imp}, for all $n$ large enough, we obtain the inequality
$$
|P_n(\beta,u, x,y)-f(y|x)|
\leq C_0(\sharp M +1)\delta n^{-\beta}\frac{|U_n(\beta,u,x)|}{\p(X\in U_n(\beta,u,x))}
$$
for $u\in (0,\delta]$, $\mu$-almost all $x\in \mathbb{R}^d$ and $y\in M$.
In view of \eqref{e_imp}, for $\mu$-almost all $x\in B_{1,n}$, $u\in (0,\delta]$ and $y\in M$,
$$
\frac{\p(X\in U_n(\beta,u,x))}{|U_n(\beta,u,x)|} \geq f_X(x) - C_0(\sharp M)u n^{-\delta}
\geq \theta_n - C_0(\sharp M)\delta n^{-\beta}\geq \frac{1}{2} n^{-\theta}
$$
if $n$ is large enough ($n\geq N(\sharp M, C_0,\delta,\beta,\theta)$)
and $\beta > \theta$. Thus, for such $n$
and indicated  $u$, $x$ and $y$,
\begin{equation}\label{est_pn}
|P_n(\beta,u, x,y)-f(y|x)|\leq 2C_0 \delta(\sharp M +1) n^{-\beta + \theta}.
\end{equation}
Note now that, for $k_n>1$ and $P_n:=P_n(\beta,u, x,y)$,
$$
\e \log (\mu_n+1)= \log((k_n -1)P_n) + \e \log\left(\frac{\mu_n+1}{(k_n -1)P_n}\right).
$$
Set $G_n:=(0,\delta]\times B_{1,n}\times M$.
We will demonstrate that
\begin{equation}\label{important}
\sup_{(u,x,y)\in G_n}\left|\e \log\left(\frac{\mu_n+1}{(k_n -1)P_n}\right)\right|\to 0,\;\;n\to \infty.
\end{equation}
According to the Lyapunov inequality it is sufficient to prove that
\begin{equation}\label{lya}
\sup_{(u,x,y)\in G_n}\e \left(\log\left(\frac{\mu_n+1}{(k_n -1)P_n}\right)\right)^2\to 0,\;\;n\to \infty.
\end{equation}
Introduce
$
\eta_n:= \frac{\mu_n-(k_n-1)P_n +1}{(k_n-1)P_n}.
$
Then
$$
\e \left(\log\left(\frac{\mu_n +1}{(k_n-1)P_n}\right)\right)^2 = \e \left( \log \left(1 + \eta_n\right)\right)^2
$$
$$
=\e \left( \left(\log (1 + \eta_n\right))^2\ind\left\{|\eta_n|<\frac{1}{2}\right\}\right)
+
\e \left( \left(\log (1 + \eta_n\right))^2\ind\left\{|\eta_n|\geq\frac{1}{2}\right\}\right):= T_1(n)+T_2(n)
$$
where $T_1(n)=T_1(n;k_n,u,x,y,\alpha)$, $T_2(n)=T_2(n;k_n,u,x,y,\alpha)$.
For $k_n>1$ and $(u,x,y)\in G_n$,
$
\frac{1}{(k_n-1)P_n}\leq \frac{\mu_n+1}{(k_n-1)P_n}\leq \frac{2}{P_n}.
$
Consequently,
$$
\left|\log\left(\frac{\mu_n+1}{(k_n-1)P_n}\right)\right|\leq
\max\left\{
\left|\log \left(\frac{1}{(k_n-1)P_n}\right)\right|, \left|\log \left(\frac{2}{P_n}\right)\right|\right\}.
$$
Taking into account \eqref{est_pn} and the bound $f(y|x)> n^{-\nu}$
for $(u,x,y)\in G_n$, we see that if
$0<\nu<\beta-\theta$
then $ \frac{1}{2}n^{-\nu}\leq P_n\leq 1$ and $\frac{1}{4} n^{- \nu}k_n\leq (k_n-1)P_n\leq k_n$ when $n$ is large enough.
Therefore, for all $n$ large enough,
$$
\left|\log\left(\frac{\mu_n+1}{(k_n-1)P_n}\right)\right|\leq b \log n
$$
where $b=b(\alpha,\nu)$ does not depend on $n$.
If $0<\nu<\alpha$ then, for all $n$ large enough,
$$
T_2(n)\leq (b \log n)^2 \p\left(|\eta_n|\geq \frac{1}{2}\right)\leq (b \log n)^2 \p(|\mu_n - (k_n-1)P_n|\geq \frac{1}{4}(k_n-1)P_n).
$$
For a random variable $H(m)\sim Bin(m,p)$, $p \in (0,1)$, $p=p_m$ and $\varepsilon=\varepsilon_m >0$, the Hoeffding inequality (see, e.g., \cite{Massart}, p. 22, and further generalizations
there) yields
$$
\p(|H(m) - mp| \geq m\varepsilon) \leq 2 \exp\{-2m\varepsilon^2\}.
$$
We employ this inequality for $m=(k_n-1)$, $k_n\propto n^{\alpha}$, $p=P_n(\beta,u, x,y)$
and  $ \varepsilon_{k-1}=\frac{P_n}{4}$. Then
$$
\sup_{(u,x,y)\in G_n}T_2(n)\leq 2(b \log n)^2 \exp\left\{-\frac{1}{8}(k_n-1)P_n^2\right\}\to 0,\;\;n\to \infty,
$$
whenever $\alpha - 2\nu>0$ (we can take positive $\nu$ arbitrary small).

To get an upper bound for $T_1(n)$ we note that $|\log(1+ z)|\leq 2|z|$ for $|z|<\frac{1}{2}$.
Hence
$$
T_1(n) = \e \left( \left(\log (1 + \eta_n\right))^2\ind\left\{|\eta_n|<\frac{1}{2}\right\}\right)
\leq 4 \e \eta_n^2.
$$
It holds
\begin{equation}\label{eq_7}
\e\eta_n^2 = \frac{(k_n-1)P_n(1-P_n) +1}{((k_n-1)P_n)^2}\leq \frac{1}{(k_n-1)P_n}
+\frac{1}{((k_n-1)P_n)^2}.
\end{equation}
We have seen that, for $(u,x,y)\in G_n$ and all $n$ large enough, the following inequality takes place
 $(k_n-1)P_n \geq \frac{1}{4}n^{- \nu}k_n\to \infty$ if $0<\nu<\alpha/2$
(we also assume that $\nu <\beta - \theta$). Therefore, the right-hand side of \eqref{eq_7}
tends to zero as $n\to \infty$.
Thus we have verified that
$
\sup_{(u,x,y)\in G_n}T_1(n)\to 0$, $n\to \infty$.
In such a way \eqref{lya} and \eqref{important} are proved.
Hence
	\begin{gather*}
\sup_{(u,x,y)\in G_n}	|\e\log(\mu_{n}+1) - \log((k_n-1)P_n(\beta,u,x,y))| \to 0,\;\;n\to \infty.
	\end{gather*}
	
Introduce notation $F_n := \frac{P_n(\beta, u, x, y)}{f(y|x)}$ where $f(y|x)>0$.
Then
	\begin{gather*}
	\sup_{(u,x,y)\in G_n}	\left|F_n - 1\right| = \sup_{(u,x,y)\in G_n}	\frac{|P_n(\beta,u,x,y) - f(y|x)|}{f(y|x)}	
	\leq cn^{-\beta+\theta+\nu} \to 0,\;\;n\to \infty,
	\end{gather*}
if $0<\nu<\beta-\theta$ and $c$ is defined by means of
\eqref{est_pn} and does not depend on $n$.
We see that $\sup_{(u,x,y)\in G_n}	\left|F_n - 1\right|<\frac{1}{2}$ for all
$n$ large enough.
Then
$$
	|\log ((k_n-1)P_n(\beta,u,x, y)) - \log ((k_n-1)f(y|x))|  = |\log (1+(F_n-1))|\leq 2|F_n-1|.
$$
So we come to the relation $
	\sup_{(u,x,y)\in G_n}|\log F_n| \to 0,\;\;n\to \infty.
	$
Thus
\begin{equation}\label{MU}
	\sup_{(u,x,y)\in G_n}|\e\log(\mu_{n}+1) - \log((k_n-1)f(y|x)| \to 0,\;\;n\to \infty.
\end{equation}
In a similar way we verify that
$
	\sup_{(u,x,y)\in G_n}|\e\log(\mu_{n}+2) - \log((k_n-1)f(y|x)| \to 0,\;\;n\to \infty
$.
Taking into account \eqref{st_2} we ascertain that
	\begin{multline*}
	\sup_{(u,x,y)\in G_n} |\e \left( \log  (\xi_{n,k,1}(x,y)+1)  |\rho_{n,k,1}(x) = un^{-\beta} \right) -  \log((k-1)f(y|x)| \\
	\leq \sup_{(u,x,y)\in G_n}|\e\log(\mu_{n}+1) - \log((k-1)f(y|x)| + \sup_{(u,x,y)\in G_n}|\e\log(\mu_{n}+2) - \log((k-1)f(y|x)| \to 0
	\end{multline*}
as $n\to \infty$.
Consequently,
	\begin{gather*}
	\sum_{y\in M} \int_{B_{1,n}} \int_{(0,\delta]} \left( \e \left( \log\left( \frac{\xi_{n,k,1}(x,y)+1}{k}\right)   \Big|\rho_{n,k,1}(x) = un^{-\beta} \right) -  \log\left(\frac{k-1}{k}\right) -
	\log f(y|x) \right)\\
\times  f^{(k)}_{n,x,\beta}(u)du
	f(x,y)dx \\ =
	S_2(n,k) - \sum_{y\in M} \int_{B_{1,n}} \int_{(0,\delta]} \left(\log\left(\frac{k-1}{k}\right) +
	\log f(y|x) \right) f^{(k)}_{n,x,\beta}(u)du f(x,y)dx\to 0
	\end{gather*}
as $n\to \infty$ (recall that $k=k_n$).	It remains to show that
	$$
	-\sum_{y\in M} \int_{B_{1,n}} \int_{(0,\delta]} \left( \log\left(\frac{k-1}{k}\right) +
	\log f(y|x) \right) f^{(k)}_{n,x,\beta}(u)du f(x,y)dx \to H(Y|X),\;\;n\to \infty.
	$$
Firstly, we can write
	\begin{gather*}
	0 \leq -\sum_{y\in M} \int_{B_{1,n}} \int_{(0,\delta]} \log\left(\frac{k-1}{k}\right) f^{(k)}_{n,x,\beta}(u)du f(x,y)dx \\
	= -\log\left(1 - \frac{1}{k}\right) \sum_{y\in M} \int_{B_{1,n}} \int_{(0,\delta]} f^{(k)}_{n,x,\beta}(u)du f(x,y)dx \leq -\log\left(1 - \frac{1}{k}\right) \to 0,\;\;n\to \infty.
	\end{gather*}
Secondly,
	\begin{gather*}
	\sum_{y\in M} \int_{B_{1,n}} \int_{(0,\delta]} \log f(y|x) f^{(k)}_{n,x,\beta}(u)du f(x,y)dx = \sum_{y\in M} \int_{B_{1,n}} \log f(y|x)  \int_{(0,\delta]} f^{(k)}_{n,x,\beta}(u)du f(x,y)dx.
	\end{gather*}
Inequality \eqref{ineq4} yields
$$
0\leq \Delta_n:=	\sup_{x \in B_{1,n}} \left(1-\int_{(0,\delta]} f^{(k)}_{n,x,\beta}(u)du\right) \to 0,\;\;n\to \infty.
	$$
Thus
$$
\left|\sum_{y\in M} \int_{B_{1,n}} \int_{(0,\delta]} \log f(y|x) f^{(k)}_{n,x,\beta}(u)du f(x,y)dx
- \sum_{y\in M} \int_{B_{1,n}}  \log f(y|x)  f(x,y)dx\right|
$$
$$
\leq \Delta_n \sum_{y\in M} \int_{\mathbb{R}^d} |\log f(y|x)|  f(x,y)dx \to 0,\;\;n\to \infty.
$$
Note now that
$$
\sum_{y\in M} \int_{B_{1,n}}\log f(y|x) f(x,y)dx \to \sum_{y\in M} \int_{\mathbb{R}^d}\log f(y|x) f(x,y)dx,\;\;n\to \infty,
$$
since $B_{1,n} \nearrow \{(x,y): f_X(x)>0, \cap_{y\in M}\{f(y|x)>0\}\}$ and $\int_B h(x)dx=0$
when $h\in L^1(\mathbb{R})$ and $B$ is a Borel subset of $\mathbb{R}^d$ such that $\mu(B)=0$.
Consequently,
	$$
	-\sum_{y\in M} \int_{B_{1,n}} \int_0^\delta \log f(y|x) f^{(k)}_{n,x,\beta}(u)du f(x,y)dx \to H(Y|X),\;\;n\to \infty.
	$$
To prove Theorem \ref{th_1} we have imposed on parameters
$\beta>0$, $\theta>0$,
and $\nu>0$ the following conditions:
$\beta >\theta$, $\nu<\beta - \theta$, $(d+1)\theta<1-\alpha$,
$\alpha - 2\nu>0$.
For each given $\alpha\in (0,1)$ we can guarantee the validity of the indicated inequalities.
Namely, one can pick
$\beta \in (0,\frac{1-\alpha}{d+1})$ and then take $\theta \in (0,\beta)$.
After that it remains to fix $\nu \in (0,\beta-\theta)$ so that $\nu < \frac{1}{2}\alpha$.

Thus the proof of Theorem \ref{th_1} is complete. $\square$
	
\vskip0.2cm
\noindent

\vskip0.3cm
{\large {\bf 4. Proof of Theorem \ref{th_2}}}	
\vskip0.2cm	
	
The proof is divided into several steps. Since $Z_1,Z_2,\ldots $ are i.i.d. observations, one has
\begin{gather*}
\e (\ent - H(Y|X))^2 = \frac{1}{n} \e (\widehat{H}_{n,k,1} - H(Y|X))^2\\ + \left(1 - \frac{1}{n}\right) \e (\widehat{H}_{n,k,1} - H(Y|X))(\widehat{H}_{n,k,2} - H(Y|X)).
\end{gather*}
We will see that the expectations in the right-hand side of the latter formula are finite. Moreover, we will verify that, for $n \to \infty$,
\vskip0.2cm	
(A)
\;\; $\e (\widehat{H}_{n,k,1} - H(Y|X))(\widehat{H}_{n,k,2} - H(Y|X)) = o(1)$,
\vskip0.2cm	
(B)\;\; $\e(\widehat{H}_{n,k,1} - H(Y|X))^2 = o(n)$.
\vskip0.4cm	

\noindent
Consider $n>2$, $k\in \{1,\ldots,n\}$, $z_j = (x_j, y_j) \in \mathbb{R}^d \times M$, $j \in \{1,2\}$. Set
$$
\widehat{H}_{n,k,j}(z_1,z_2)
:= -\log  \left(\frac{\xi_{n,k,j}(z_1,z_2, Z_3, \dots, Z_n)+1}{k}\right).
$$
The independence of observations $Z_1,\ldots,Z_n$ implies that
\begin{gather*}
\e \left((\widehat{H}_{n,k,1} - H(Y|X))(\widehat{H}_{n,k,2} - H(Y|X))|Z_1=z_1, Z_2=z_2\right) \\
= \e
(\widehat{H}_{n,k,1}(z_1, z_2) - H(Y|X))(\widehat{H}_{n,k,2}(z_1,z_2) - H(Y|X))=:\mathcal{H}_{n,k}(z_1,z_2).
\end{gather*}
Therefore,
\begin{gather*}
\e (\widehat{H}_{n,k,1} - H(Y|X))(\widehat{H}_{n,k,2} - H(Y|X))  \\
=\sum_{y_1, y_2 = 1}^{m} \int\limits_{\mathbb{R}^d} \int\limits_{\mathbb{R}^d}
\mathcal{H}_{n,k}(z_1,z_2) f(x_1,y_1) f(x_2,y_2) \; dx_1 dx_2.
\end{gather*}

Due to the de la Vall\'ee Poussin theorem (see, e.g., \cite{Borkar}, p. 10), for establishing (A) it suffices to prove validity of
the following two statements.

\noindent
1) If $dQ(x_1,x_2):= f(x_1,y_1) f(x_2,y_2) \; dx_1 dx_2$, then for
each $y_1, y_2 \in M$ and $Q$-a.s. $(x_1, x_2) \in \mathbb{R}^d\times \mathbb{R}^d$,

\begin{equation}\label{B.1}
	\mathcal{H}_{n,k}(z_1,z_2)
	\to h(x_1,y_1)h(x_2,y_2),\;\;n \to \infty,
\end{equation}
here $h(x,y):= -\log f(y|x)-H(Y|X)$, $x\in \mathbb{R}^d$, $y\in M$, $z_j = (x_j, y_j), j \in \{1, 2\}$.

\noindent
2) For some $a > 0$,
	\begin{equation}\label{B.2}
	\sup_n\sum_{y_1, y_2 \in M} \int\limits_{\mathbb{R}^d} \int\limits_{\mathbb{R}^d} \left|\mathcal{H}_{n,k}(z_1,z_2)  \right|^{1+a} f(x_1,y_1) f(x_2,y_2) \, dx_1 dx_2 < \infty.
	\end{equation}
Indeed, $\eqref{B.1}$ and $\eqref{B.2}$ imply that
$$
\sum_{y_1, y_2 \in M} \int\limits_{\mathbb{R}^d} \int\limits_{\mathbb{R}^d} \mathcal{H}_{n,k}(z_1,z_2)   f(x_1,y_1) f(x_2,y_2)\; dx_1 dx_2
$$
$$
\to \sum_{y_1, y_2 \in M} \int\limits_{\mathbb{R}^d} h(x_1,y_1)f(x_1,y_1)\,dx_1 \int\limits_{\mathbb{R}^d} h(x_2,y_2) f(x_2,y_2)\,dx_2 =0,\;\;n\to \infty.
$$
In view of the Jensen conditional inequality it is easily seen that \eqref{B.2} holds if
\begin{equation*}
\sup_n\e \left|(\widehat{H}_{n,k,1} - H(Y|X))(\widehat{H}_{n,k,2} - H(Y|X))\right|^{1+a} < \infty.
\end{equation*}
The Cauchy -  Schwartz inequality yields
\begin{gather*}
\e \left|(\widehat{H}_{n,k,1} - H(Y|X))(\widehat{H}_{n,k,2} - H(Y|X))\right|^{1+a}
\leq
\e |\widehat{H}_{n,k,1} - H(Y|X)|^{2+2a}.
\end{gather*}
Thus, as $a>0$ can be taken arbitrary small, \eqref{B.2} holds if, for some $\varepsilon >0$,
\begin{equation}
\sup_n \e |\widehat{H}_{n,k,1} - H(Y|X)|^{2+\varepsilon} < \infty. \label{sup_bound}
\end{equation}
On applying  the Lyapunov moment inequality we observe that \eqref{sup_bound} guarantees validity of (B) .

Employing the reasoning used to prove Theorem \ref{th_1} one arrives at an expression
\begin{multline}
\e \left|-\log \left(\frac{\xi_{n,k,1}(Z_1, \dots, Z_n)+1}{k}\right) - H(Y|X)\right|^{2+\varepsilon} \\
= \sum_{y \in M} \int\limits_{\mathbb{R}^d} \e \left|\log  \left(\frac{\xi_{n,k,1}(x, y)+1}{k}\right) + H(Y|X) \right|^{2+\varepsilon} f(x,y) \, dx \label{initial_expectation}.
\end{multline}
Fix an arbitrary $\beta >0$ and, for $x\in \mathbb{R}^d$, $y\in M$ and $u>0$, set
$$
I_{n,k,\beta,\varepsilon}(x, y, u):=
\e \left( \left. \left| \log  \left(\frac{\xi_{n,k,1}(x,y)+1}{k}\right) + H(Y|X) \right|^{2+\varepsilon} \right|n^{\beta}\rho_{n,k,1} = u \right).
$$
Then we can write
$$
\e \left|\log  \left(\frac{\xi_{n,k,1}(x,y)+1}{k}\right) + H(Y|X) \right|^{2+\varepsilon}
=
\int_{0}^{\infty} I_{n,k,\beta,\varepsilon}(x, y, u) f^{(k)}_{n,x,\beta}(u)\, du
$$
where $f^{(k)}_{n,x,\beta}(\cdot)$ is a density of random variable $n^{\beta}\rho_{n,k,1}$
(see formula \eqref{dens} in Appendix).
Taking into account the finiteness of the set $M$ it is sufficient to prove that, for each $y \in M$, one has $\sup_n |\mathcal{I}_{n,k,i}(y)| < \infty$ where
\begin{gather*}
\mathcal{I}_{n,k,i}(y):= \int\limits_{B_{i, n}} \int_{0}^{\infty} I_{n,k,\beta,\varepsilon}(x, y, u)  f^{(k)}_{n,x,\beta}(u) f(x,y) \, du dx, \;\;\; i =1,2,
\end{gather*}
and $B_{i, n}$ were defined in \eqref{B1_B2}.
Clearly, if $k>1$ then, for any  $x\in \mathbb{R}^d$ and $y\in M$, the following is $\p$-a.s. true
\begin{equation}\label{b0}
\left|\log \left( \frac{\xi_{n,k,1}(x,y)+1}{k}\right) + H(Y|X) \right|\leq \log k + |H(Y|X)|.
\end{equation}
Hence,
using \eqref{dens_1} and \eqref{dens_3}
with $2+\varepsilon$ instead of $1+\varepsilon$,
we come to the relation
 $\sup_n |\mathcal{I}_{n,k,2}| < \infty$.

Fix $\delta >0$. One has
$
\mathcal{I}_{n,k,1}(y) = S_{n,k,1}(y) + S_{n,k,2}(y)$ where, for $V_1=(\delta,\infty)$ and
$V_2=(0,\delta]$,
$$
S_{n,k,j}(y) = \int\limits_{B_{1,n}} \int_{V_j} I_{n,k,\beta,\varepsilon}(x, y, u) f^{(k)}_{n,x,\beta}(u) f(x, y) \, du dx,\;\;j=1,2.
$$
In similarity  to  \eqref{ineq1}, for each $y\in M$ and all $n$ large enough,
basing on \eqref{b0} we obtain that
\begin{gather*}
0 \leq S_{n,k,1}(y) \leq (\log k+|H(Y|X)|)^{2+\varepsilon}  \int\limits_{B_{1,n}}
\int_\delta^\infty f^{(k)}_{n,x,\beta}(u)f(x,y) \, du \, dx\\
\leq (\log k+|H(Y|X)|)^{2+\varepsilon}  \int\limits_{B_{1,n}} \p(\eta_n(\delta,\beta, x) \leq k-1) f_X(x)\, dx.
\end{gather*}
According to inequality \eqref{ineq4} one has
\begin{equation}\label{ineq6}
0 \leq S_{n,k,1}(y) \leq \frac{(\log k+|H(Y|X)|)^{2+\varepsilon}}{\sqrt{2 \pi R_n}} e^{-R_n/2} \to 0,\;\;n\to \infty,
\end{equation}
where $R_n$ was introduced in \eqref{R_cond}.

Now we turn to the estimation of $S_{n,k,2}(y)$. Lemma \ref{lem_1} yields that
\begin{multline*}
I_{n,k,\beta,\varepsilon}(x, y, u)
= \e \left| \log  \left(\frac{\mu_n+1}{k}\right)  + H(Y|X) \right|^{2+\varepsilon} \p(Y \ne y | \|X - x\| = un^{-\beta}) \\
+ \e \left| \log  \left(\frac{\mu_n+2}{k}\right) + H(Y|X) \right|^{2+\varepsilon}  \p(Y = y | \|X - x\| = un^{-\beta}),
\end{multline*}
here the variable $\mu_n$ is defined in \eqref{mu_n} and does not depend on $(X,Y)$. We show that
\begin{gather}\label{ineq5}
\sup_{(u, x, y) \in G_n} \left|I_{n,k,\beta,\varepsilon}(x, y, u) - |\log f(y|x) + H(Y|X)|^{2+\varepsilon} \right| \to 0, \;\;n\to \infty,
\end{gather}
where $G_n= (0,\delta]\times B_{1,n}\times M$.
It is sufficient to prove that, as $n\to \infty$,
\begin{equation}\label{L2_ineq1}
\sup_{(u, x, y) \in G_n}\left| \e \left| \log   \left(\frac{\mu_n+1}{k}\right)  + H(Y|X) \right|^{2+\varepsilon} - |\log f(y|x) + H(Y|X)|^{2+\varepsilon} \right| \to 0,
\end{equation}
\begin{equation}\label{L2_ineq2}
\sup_{(u, x, y) \in G_n}\left| \e \left| \log  \left(\frac{\mu_n+2}{k}\right)  + H(Y|X) \right|^{2+\varepsilon} - |\log f(y|x) + H(Y|X)|^{2+\varepsilon} \right| \to 0.
\end{equation}
For $\nu \in (0,1/2)$ and $\gamma \in (0, 1/2 - \nu)$
one has $\e \left| \log  \frac{\mu_n+1}{k} + H(Y|X) \right|^{2+\varepsilon}	=  T_{n,1} + T_{n,2}$ where
\begin{gather*}
T_{n,1} :=  \e \left(\left| \log  \left(\frac{\mu_n+1}{k}\right) + H(Y|X) \right|^{2+\varepsilon}  \mathbb{I}\{ |\mu_n - \e \mu_n| \leq (k-1)^{1/2 + \gamma}  \}\right), \\
T_{n,2} :=  \e \left(\left| \log  \left(\frac{\mu_n+1}{k}\right) + H(Y|X) \right|^{2+\varepsilon} \mathbb{I}\{ |\mu_n - \e \mu_n| > (k-1)^{1/2 + \gamma}  \}\right),
\end{gather*}
here $T_{n,j}=T_{n,j}(\beta,\varepsilon, \gamma, u,x,y)$, $j=1,2$.
Note that by elementary properties of the binomial distribution
$
\e \mu_n = (k-1)P_n(\beta, u, x, y)$ with $P_n(\beta, u, x, y)$ introduced in
\eqref{Pn}.

In view of \eqref{b0}, for $(u,x,y)\in G_n$ and all $n$ large enough, one can write
$$
0 \leq \left| \log \left(\frac{\mu_n+1}{k}\right) + H(Y|X) \right|^{2+\varepsilon} \leq (\log k+|H(Y|X)|)^{2+\varepsilon}.
$$
Applying the Hoeffding inequality with $m = k-1$, $H(m) = \mu_n$, $p = P_n$, $\varepsilon = (k-1)^{-1/2 + \gamma}$ we have
$$
0 \leq T_{n, 2} \leq (\log k +|H(Y|X)|)^{2+\varepsilon}  \p(|\mu_n - \e \mu_n| > (k-1)^{1/2+\gamma} )
$$
\begin{equation}\label{T_2}
\leq (\log k+|H(Y|X)|)^{2+\varepsilon}  2e^{-2(k-1)^{2\gamma}} \to 0,\;\;n\to \infty.
\end{equation}
Note also that
\begin{equation}\label{ABC}
|\log f(y|x) + H(Y|X)|^{2+\varepsilon} \p(|\mu_n - \e \mu_n| > (k-1)^{1/2+\gamma} )\to 0,\;\;n\to
\infty.
\end{equation}

Set $B_i(n,p) = \p(\mu_n = i)$, $i=0,\ldots,k-1$, and write $T_{n,1}$ in the following way
\begin{equation}\label{T1}
T_{n,1} = \sum_{i:|i - \e \mu_n| \leq (k-1)^{1/2 + \gamma}} \left| \log  \left(\frac{i+1}{k}\right) + H(Y|X) \right|^{2+\varepsilon} B_i(n,p).
\end{equation}
To get the upper bound for $|T_{n,1} - |\log f(y|x) + H(Y|X)|^{2+\varepsilon}|$ we employ the Lagrange formula for a function
$
g(z) = |\log z + a|^{2+\varepsilon}$, $z > 0$, $a \in \mathbb{R}$.
For $z,z_0>0$, one has
$$
g(z)-g(z_0)= g'(\xi)(z-z_0),\;\;\xi= z_0+\lambda(z-z_0),\;\;\lambda=\lambda(z,z_0) \in (0,1).
$$
Take $z = \frac{i+1}{k}$, $z_0 = f(y|x)$ and $a = H(Y|X)$. Then, for $i$ belonging to the summation set in \eqref{T1}, in view of \eqref{est_pn} we get
$$
|z - z_0| =
\left| \frac{i+1}{k} - f(y|x) \right| \leq
\left| \frac{i - (k-1)P_n}{k} \right| + \frac{1}{k}
+ \left|\frac{(k-1)P_n}{k} - P_n \right|
+ \left| P_n - f(y|x) \right|
$$
$$
\leq (k-1)^{-1/2+\gamma} + 2k^{-1} + 2C \delta(\sharp M +1) n^{-\beta + \theta} := \mathcal{Z}_n \propto n^{-\beta+\theta},
$$
provided that $\beta - \theta < 1/2-\gamma$.
Note that $
g'(z) = (2+\varepsilon)z^{-1}(\log z + a)|\log z + a|^{\varepsilon},\;\;z>0.$
Hence $|g'(\xi)|\leq (2+\varepsilon)|\xi|^{-1}|\log \xi +a|^{1+\varepsilon}$.
For all $n$ large enough,
$$
|\log \xi|\leq \max\{|\log z_0|,|\log z|\} \leq \max\{|\log \nu_n|, \log k\}
\leq c\log n
$$
where $c=c(\alpha,\nu)$.
Clearly, $|\xi|\geq \min\{z,z_0\}$. In our case
$z_0 \geq n^{-\nu}$ and $z \geq P_n(k_n-1)-(k_n-1)^{\frac{1}{2}+\gamma}$.
According to \eqref{est_pn} we can write
$P_n \geq \frac{1}{2}n^{-\nu}$ for all $n$ large enough
if $\nu<\beta -\theta$. Therefore, for all $n$ large enough, one can see that
$z \geq c_1 n^{-\nu+\alpha}$ if $\alpha - \nu > \frac{1}{2}+\gamma$ (here $c_1=c_1(\alpha)$).
The latter inequality holds if $\nu<2\alpha$ (then we take positive $\gamma$ which is small enough).  Thus $|\xi|^{-1}
\leq n^{\nu}$ for all $n$ large enough.
Consequently, uniformly in $i$ belonging to the summation set in definition of $T_{n,1}$,
$$
|g(z) - g(z_0)| = |g'(\xi)||z-z_0|\leq (2+\varepsilon)(ñ\log n)^{1+\varepsilon}n^{\nu}\mathcal{Z}_n \propto
(\log n)^{1+\varepsilon} n^{\nu-\beta +\theta}\to 0,\;\;n\to \infty,
$$
whenever $\nu<\beta-\theta$.
Taking into account that $\sum_{i:|i-\e\mu_n|\leq (k-1)^{1/2+\gamma}}B_i(n,p)\leq 1$
we come to relation
\begin{gather*}
\sup_{(u,x,y)\in G_n}|T_{n,1} - |\log f(y|x) + H(Y|X)|^{2+\varepsilon}| \to 0,\;\;n\to \infty.
\end{gather*}

This formula, \eqref{T_2} and \eqref{ABC} yield \eqref{L2_ineq1}.
Relation \eqref{L2_ineq2}
is proved analogously. In such a way we establish \eqref{ineq5}. Hence,
\begin{gather*}
S_{n,k,2}(y) - \int_{B_{1, n}}\int_0^\delta |\log f(y|x) + H(Y|X)|^{2+\varepsilon} f^{(k)}_{n,x,\beta}(u) f(x, y) \, du dx \to 0,\;\;n\to \infty.
\end{gather*}
Obviously, for each $n\in \mathbb{N}$ and any $y\in M$, we get
\begin{gather*}
\int_{B_{1,n}}\int_0^\delta |\log f(y|x) + H(Y|X)|^{2+\varepsilon} f^{(k)}_{n,x,\beta}(u) f(x, y) \, du dx \\
\leq \int_{\mathbb{R}^d \cap \{x: f(y|x)>0\}} |\log f(y|x) + H(Y|X)|^{2+\varepsilon} f(x, y) \, dx < \infty.
\end{gather*}
Indeed, in view of the Minkowski inequality and since $f(y|x) = \p(Y=y|X=x) \leq 1$ for each $y \in M$ and $\p_X$-almost all $x$, it is enough to show that
$$
\int_{\mathbb{R}^d \cap \{x: f(y|x)>0\}} \left(\log \frac{1}{f(y|x)}\right)^{2+\varepsilon} f(x, y) \, dx < \infty,\;\;y\in M.
$$
For any $\varepsilon >0$, there exists such $\mathcal{T}=\mathcal{T(\varepsilon)} > 1$ that $(\log t)^{2+\varepsilon} \leq t$ whenever $t > \mathcal{T}$. Hence, for each $y\in M$,  the latter integral can be written as follows
\begin{gather*}
\int_{\{x: f(y|x) \geq 1/\mathcal{T}\}} \!\!\left(\log \frac{1}{f(y|x)}\right)^{2+\varepsilon} \!\!f(x,y) \, dx
\!+\! \int_{\{x: 0<f(y|x) < 1/\mathcal{T}\!\}}\!\! \left(\log \frac{1}{f(y|x)}\right)^{2+\varepsilon}\!\! f(y|x) f_X(x) \, dx \\
\leq (\log \mathcal{T})^{2+\varepsilon}\p(Y=y) + \int_{\{x: f(y|x) < 1/T\}} (f(y|x))^{-1} f(y|x) f_X(x) \, dx \leq (\log \mathcal{T})^{2+\varepsilon} + 1 < \infty.
\end{gather*}

Consequently, for $y \in M$,
one has
$\sup_n S_{n,k,2}(y) < \infty$. Thus by applying \eqref{ineq6} we ensure, for each $y\in M$, that
$
\sup_n |\mathcal{I}_{n,k,1}(y)| < \infty.
$
Relation $\eqref{sup_bound}$ is established.

Now we concentrate on the proof of \eqref{B.1}.
For $z_1,z_2\in \mathbb{R}^d \times M$, $n>2$ and $k=k_n$, introduce
$$
R_{n,k}(z_1,z_2):= (\widehat{H}_{n,k,1}(z_1, z_2) - H(Y|X))(\widehat{H}_{n,k,2}(z_1,z_2) - H(Y|X)).
$$
Write $\rho_{n,k,j}:=\rho_{n,k,j}(x_1,x_2)$, $j=1,2$. Fix $\delta > 0$ and set
\begin{gather*}
\mathcal{U}_{n,k,1}(z_1,z_2) := \e(
R_{n,k}(z_1,z_2)
|\rho_{n,k,1}n^{\beta} \leq \delta, \rho_{n,k,2}n^{\beta} \leq \delta) \p(\rho_{n,k,1}n^{\beta} \leq \delta, \rho_{n,k,2}n^{\beta} \leq \delta), \\
\mathcal{U}_{n,k,2}(z_1,z_2) := \mathcal{H}_{n,k}(z_1,z_2) - \mathcal{U}_{n,k,1}(z_1,z_2)
\end{gather*}
where, for an integrable random variable $\xi$, one has    $\e(\xi|A):=0$ if $\p(A)=0$.
For all $n\geq N_1$, where $N_1=N_1(H(Y|X))$, and any $z_j\in \mathbb{R}^d\times M$,  $j=1,2$,
the inequality $|R_{n,k}(z_1,z_2)|\leq 2(\log k)^2$ holds with probability one.
Therefore
$$
|\mathcal{U}_{n,k,2}(z_1,z_2)| \leq 2(\log k)^2  \p\left( \left\{ \rho_{n,k,1} \!>\!
\delta n^{-\beta} \right\} \cup \left\{ \rho_{n,k,2} \!>\! \delta n^{-\beta} \right\} \right)
%\\
\leq 4(\log k )^2  \p\left( \rho_{n,k,1} \!>\! \delta n^{-\beta} \right).
$$
Due to \eqref{f1} and \eqref{ineq4} we know, for each $z_j\in \mathbb{R}^d\times M$,  $j=1,2$, that $
|\mathcal{U}_{n,k,2}(z_1,z_2)| \to 0,\;\;n\to \infty.
$
It remains to show that, for any $y_1, y_2 \in M$ and $Q$-almost all $(x_1, x_2) \in \mathbb{R}^d\times \mathbb{R}^d$,
\begin{equation} \label{un_ineq}
\mathcal{U}_{n,k,1}(z_1,z_2) \to (-\log f(y_1|x_1) - H(Y|X))(-\log f(y_2|x_2) - H(Y|X)),\;\;n\to \infty.
\end{equation}

In the proof of Lemma \ref{lem_2} (see the Appendix) it is shown that, for each $x_1,x_2\in \mathbb{R}^d$ ($x_1\neq x_2$), a finite measure
$
G(B) := \p((\rho_{n,k,1}, \rho_{n,k,2}) \in B),
$ defined for
$B \in \mathcal{B}(\mathbb{R}^2) \cap D_{x_1, x_2}$, where
$D_{x_1, x_2} = (0,|x_1-x_2|/3]\times (0,|x_1-x_2|/3]$ can be written in the following way
$$
\p((\rho_{n,k,1}, \rho_{n,k,2}) \in B) = \int_B g_{n,k}(x_1, x_2, u_1, u_2) \, du_1 du_2.
$$
Here $g_{n,k}(x_1, x_2, \cdot, \cdot)\colon D_{x_1, x_2} \to \mathbb{R}_+$ is a certain
integrable function (w.r.t. the restriction to $D_{x_1,x_2}$ of $mes \otimes mes$).
 Clearly, for $0<k<n-2$,
\begin{gather*}
\p(\rho_{n,k,i}(x_1, x_2) = 0) =\sum_{j=k}^{n-2} \binom{n-2}{j} \p( \|X - x_i \| = 0)^j (1 - \p( \|X - x_i \| = 0))^{n-2-j} = 0,\;\; i=1,2.
\end{gather*}
since $X$ has a density (and $k = k(n) \propto n^\alpha$, $\alpha \in (0,1)$).

Let us set
$g_{n,k}(x_1, x_2, \cdot, \cdot)=0$ on $\overline{D}_{x_1, x_2} \setminus D_{x_1, x_2}$ where $\overline{D}_{x_1, x_2}$ stands for a closure of $D_{x_1, x_2}$ in $\mathbb{R}^2$.
Take $N(x_1,x_2,\beta,\delta)$ to ensure that $\delta n^{-\beta} \leq |x_1-x_2|/6$ for all $n\geq N(x_1,x_2,\beta,\delta)$. Then $Q_{n,\beta,\delta}:=[0,\delta n^{-\beta}]\times
[0,\delta n^{-\beta}]\subset \overline{D}_{x_1, x_2}$ and
$B(x_1,u_1n^{-\beta})\cap B(x_2,u_1n^{-\beta})=\varnothing$ for $u_i\in [0,\delta]$, $i=1,2$.
Thus, for any  $B \in \mathcal{B}(\mathbb{R}^2)$ and $n \geq  N(x_1, x_2,\beta,\delta)$,
\begin{gather*}
\p(\zeta \in B|\zeta\in Q_{n,\delta,\beta}) = \frac{\p(\zeta \in B \cap Q_{n,\delta,\beta})}{\p(\zeta \in Q_{n,\delta,\beta})}
= \frac{1}{\p(\zeta \in Q_{n,\delta,\beta})} \int_{B \cap Q_{n,\delta,\beta}} g_{n,k}(x_1, x_2, u_1, u_2) \, du_1 du_2
\end{gather*}
since $B \cap Q_{n,\delta,\beta} \subset \overline{D}_{x_1, x_2}$.
Note that $\p(\zeta \in Q_{n,\delta,\beta}) > 0$ for all
$n$ large enough as
\begin{gather*}
\p( \zeta \in Q_{n,\delta,\beta}) = \sum_{\substack{
		l_1 \geq k,
		l_2 \geq k, \\
		l_1+l_2 \leq n-2
}} \binom{n}{l_1} \binom{n-l_1}{l_2} p_{n, x_1, \delta, \beta}^{l_1} p_{n, x_2, \delta, \beta}^{l_2}
(1 - p_{n, x_1, \delta, \beta} - p_{n, x_2, \delta, \beta}))^{n-l_1-l_2}.
\end{gather*}
Here $p_{n, x_i, \delta, \beta} = \p(\|X - x_i\| \leq \delta n^{-\beta})>0$ for all $n\in \mathbb{N}$ and $i=1,2$,
as $\p_X(B(x,r))>0$ for any $x\in \mathbb{R}^d$ and $r>0$ (because $f_X(z)>0$ for
$\mu$-almost all $z\in \mathbb{R}^d$).
We also take into account that $p_{n, x_i, \delta, \beta}\to 0$ as $n\to \infty$, $i=1,2$.

Therefore a function $\tilde{g}_{n,k}(x_1, x_2, \cdot, \cdot): \mathbb{R}^2 \to \mathbb{R_{+}}$,
$$
\tilde{g}_{n,k}(x_1, x_2, u_1, u_2) = \begin{cases}
\frac{1}{\p(\zeta \in Q_{n,\delta,\beta})} g_{n,k}(x_1, x_2, u_1, u_2), \,\,\, \text{if} \,\, (u_1, u_2) \in Q_{n,\delta,\beta}, \\
0, \,\,\, \text{if} \,\, (u_1, u_2) \in \mathbb{R}^2 \setminus Q_{n,\delta,\beta},
\end{cases}
$$
is a probability density of the measure $\p(\zeta \in \cdot|\zeta \in Q_{n,\delta,\beta})$ which is defined on  $(\mathbb{R}^2, \mathcal{B}(\mathbb{R}^2))$.
Thus the measure $\p(n^{\beta} \zeta \in \cdot|\zeta \in Q_{n,\delta,\beta}))$  on this space
has a density (w.r.t. a restriction of $mes \otimes mes$ on $(\mathbb{R}^2, \mathcal{B}(\mathbb{R}^2))$)
$
f_{n, k, \beta,\delta}(x_1, x_2, \cdot, \cdot) \colon \mathbb{R}^2 \to \mathbb{R_{+}}.
$

Now we are going to employ Lemma \ref{lem_3} with
$$W := R_{n,k}(z_1,z_2),\;\;
V := (n^{\beta}\rho_{n,k,1}(x_1,x_2), n^{\beta}\rho_{n,k,2}(x_1,x_2)),\;\;
B := \{(u_1, u_2) \in [0,\delta]\times[0,\delta]\}.
$$
Note that $\{V \in B\} = \{\zeta \in Q_{n,\delta,\beta}  \}$ and $\e W$ exists.
Consequently, for considered $x_1, x_2$, $\beta$, $\delta$ and $n > N(x_1, x_2,\beta,\delta)$, the following formula is valid
\begin{gather*}
\mathcal{U}_{n,k,1}(z_1,z_2) = \e(W|V \in B) \, \p(\rho_{n,k,1}n^{\beta} \leq \delta, \rho_{n,k,2}n^{\beta} \leq \delta)\\
\! =\! \int_0^\delta\!\! \int_0^\delta J_{n,k,\beta}(z_1, z_2, u_1, u_2) \, d\p_{V, B}(u_1, u_2)  \, \p(\rho_{n,k,1}n^{\beta} \leq \delta, \rho_{n,k,2}n^{\beta} \leq \delta)
\end{gather*}
where $J_{n,k,\beta}(z_1, z_2, u_1, u_2)
:= \e (R_{n,k}(z_1,z_2)| V=(u_1,u_2))$.
It was shown that the measure $\p_{V, B}$ has a density
$f_{n,k,\beta,\delta} (x_1, x_1, \cdot,\cdot)$, therefore
\begin{gather*}
\mathcal{U}_{n,k,1}(z_1,z_2)\! = \! \int_0^\delta\!\! \int_0^\delta J_{n,k,\beta}(z_1, z_2, u_1, u_2) f_{n,k,\beta,\delta} (x_1, x_1, u_1, u_2) \, du_1 du_2 \, \p(\rho_{n,k,1}n^{\beta} \leq \delta, \rho_{n,k,2}n^{\beta} \leq \delta).
\end{gather*}

Now we show that uniformly for $\p_{V,B}$-almost all $(u_1, u_2) \in (0, \delta] \times (0, \delta]$
\begin{gather}\label{ineq7}
|J_{n,k,\beta}(z_1, z_2, u_1, u_2) - (-\log f(y_1|x_1) - H(Y|X))(-\log f(y_2|x_2) - H(Y|X))| \to 0,\;\;n\to \infty.
\end{gather}
Set $\xi\! :=\!(\xi_{n,k,1},\xi_{n,k,2})$ where $\xi_{n,k,1}\!=\!\xi_{n,k,1}(z_1, z_2)$, $i\!=\!1,2$.
Due to Lemma \ref{lem_4}, for $\p_V$-almost all $(u_1, u_2)$,
\begin{gather*}
J_{n,k,\beta}(z_1, z_2, u_1, u_2) =
\sum_{r_1, r_2} \e (  R_{n,k}(z_1,z_2)
| V=(u_1,u_2), \xi = (r_1,r_2))
\p(\xi=(r_1,r_2)| V=(u_1,u_2))
\end{gather*}
Note that, for $\p_{V,\xi}$-almost all $(u_1, u_2, r_1, r_2)$,
$$
\e ( R_{n,k}(z_1,z_2)
| V=(u_1,u_2), \xi = (r_1,r_2))
\p(\xi=(r_1,r_2)| V=(u_1,u_2))
$$
\begin{equation}\label{vxi_1}
= \left( -\log \left( \frac{r_1+1}{k}\right) - H(Y|X) \right) \left( -\log \left( \frac{r_2+1}{k}\right) - H(Y|X) \right):= h(r_1,r_2)
\end{equation}
because a random variable $R_{n,k}(z_1,z_2)$
is measurable w.r.t.  $\sigma$-algebra $\sigma\{V, \xi\}$. A function
$h(r_1,r_2)$ depends also on $n,k,z_1,z_2$.

Let $O \in \mathcal{B}(\mathbb{R}_{+}^2 \times M^2)$ be the set consisting of $(u_1,u_2,r_1,r_2)$ such that \eqref{vxi_1} holds. Then $\p_{V,\xi}(O)=1$. Since $M^2$ is a finite set,
$O = \bigcup_{(r_1, r_2) \in M^2} O_{r_1, r_2} \times \{(r_1, r_2)\}$ where $O_{r_1, r_2} \in \mathcal{B}(\mathbb{R}_{+}^2)$.
Note that at least one set $O_{r_1, r_2}$ is not empty, otherwise $\p((V, \xi) \in O) \neq 1$. Consequently, $\bigcup_{(r_1, r_2) \in M^2} O_{r_1, r_2} \ne \varnothing$.
If $\bigcap_{(r_1, r_2) \in M^2} O_{r_1, r_2}\neq \varnothing$ then, for each $(u_1, u_2) \in \bigcap_{(r_1, r_2) \in M^2} O_{r_1, r_2}$,

\begin{equation}\label{vxi_2}
J_{n,k,\beta}(z_1, z_2, u_1, u_2)
= \sum_{r_1, r_2}  h(r_1,r_2)
\p(\xi=(r_1,r_2) | V=(u_1,u_2)).
\end{equation}

Define the set $\widetilde{O} = (\bigcup_{(r_1, r_2) \in M} O_{r_1, r_2} \times M^2) \setminus O$. It also can be represented as
$$
\widetilde{O} = \bigcup_{(r_1, r_2) \in M^2} \widetilde{O}_{r_1, r_2} \times \{(r_1, r_2)\}, \,\,\, \widetilde{O}_{r_1, r_2} = \bigcup_{(s_1, s_2) \in M} O_{s_1, s_2} \setminus O_{r_1, r_2}.
$$

Clearly, $\p_{V,\xi}(O) = 1$ implies that $\p_V(\bigcup_{(r_1, r_2) \in M^2} O_{r_1, r_2}) = 1$. Indeed,
\begin{gather*}
\p\left(V \in \bigcup_{(r_1, r_2) \in M} O_{r_1, r_2}\right)
\geq
\p\left(\bigcup_{(r_1, r_2) \in M} \{V \in O_{r_1, r_2}, \xi = (r_1, r_2) \}\right)
= \p((V, \xi) \in O) = 1.
\end{gather*}

If the set $\widetilde{O}$ is empty then $O_{r_1, r_2} = O_{s_1, s_2}$ for $r_1 \ne s_1, r_2 \ne s_2$, so $\bigcap_{(r_1, r_2) \in M} O_{r_1, r_2} = \bigcup_{(r_1, r_2) \in M} O_{r_1, r_2} \ne \varnothing$, thus equality holds for $\p_{V}$-almost all $(u_1, u_2)$.

Let us consider the case where $\widetilde{O} \ne \varnothing$.
Introduce
$K:=\{(r_1,r_2)\in M^2: \widetilde{O}_{r_1,r_2}\neq \varnothing\}$.
Hence $K\neq \varnothing$. If $L:=M^2\setminus K\neq \varnothing$ then,
for each $(r_1,r_2)\in L$, one has
$
O_{r_1, r_2} = \bigcup_{(s_1, s_2) \in M} O_{s_1, s_2}.
$
We have seen that $\p_V(\cup_{(s_1, s_2) \in M^2} O_{s_1, s_2})=1$.
Therefore, $\p_V(O_{r_1,r_2})=1$ for each $(r_1,r_2)\in L$.
Introduce $K_1 :=\{(r_1,r_2)\in K: \p_V(\widetilde{O}_{r_1,r_2}) = 0\}$ and $K_2 := K \setminus K_1$.
If $(r_1,r_2)\in K_1$ then $\p_V(O_{r_1, r_2}) = \p_V(\bigcup_{(s_1, s_2) \in M} O_{s_1, s_2} \setminus \widetilde{O}_{r_1, r_2}) = \p_V(\bigcup_{(s_1, s_2) \in M} O_{s_1, s_2}) - \p_V(\widetilde{O}_{r_1, r_2}) = 1$.

Now we will demonstrate that if $(r_1,r_2)\in K_2$ then
$\p(\xi = (r_1, r_2)|V = (u_1, u_2)) = 0$  for  all $(u_1, u_2) \in S_{r_1,r_2}\subset \widetilde{O}_{r_1, r_2}$, $S_{r_1, r_2} \in \mathcal{B}(\mathbb{R}_+^2)$ where $\p_V(S_{r_1,r_2})=\p_V(\widetilde{O}_{r_1, r_2}) > 0$. If the latter statement is true then
(61) is valid since we come to the trivial relation $0=0$ and, consequently, we obtain the desired formula for any $(u_1,u_2) \in A:=(\cap_{(r_1,r_2)\in L \cup K_1}O_{r_1,r_2}) \cap (\cap_{(r_1,r_2)\in K_2} (O_{r_1,r_2} \cup S_{r_1,r_2}))$ where $\p_V(A)=1$ because $\p_V(O_{r_1,r_2}) = 1$ for $(r_1, r_2) \in L \cup K_1$ and
\begin{gather*}
\p_V(O_{r_1,r_2} \cup S_{r_1,r_2}) = \p_V(O_{r_1,r_2}) + \p_V(S_{r_1,r_2}) = \p_V(O_{r_1,r_2}) + \p_V(\widetilde{O}_{r_1,r_2})\\ = \p_V(O_{r_1,r_2} \cup \widetilde{O}_{r_1,r_2})
= \p_V(\cup_{(s_1, s_2) \in M} O_{s_1, s_2}) = 1.
\end{gather*}
Here we take into account that $S_{r_1,r_2} \subset \widetilde{O}_{r_1, r_2}$ and $O_{r_1, r_2} \cap \widetilde{O}_{r_1, r_2} = \varnothing$.

For each $(r_1, r_2) \in K_2$
\begin{gather}
\p(V \in \widetilde{O}_{r_1, r_2}, \xi = (r_1, r_2)) =  \int_{\widetilde{O}_{r_1, r_2}} \p(\xi = (r_1, r_2) | V = (u_1,u_2)) \, d\p_V(u_1,u_2) = 0 \label{vxi_3}
\end{gather}
because $(\widetilde{O}_{r_1, r_2} \times \{r_1, r_2\}) \cap O = \varnothing$ and $\p_{V, \xi}(O) = 1$. Invoking that $\p(\xi = (r_1, r_2) | V = (u_1,u_2)) \geq 0$ for $\p_V$-almost all $(u_1, u_2)$, from equation \eqref{vxi_3} we infer that $\p(\xi = (r_1, r_2) | V = (u_1,u_2)) = 0$ for some set $S_{r_1, r_2} \in \mathcal{B}(\mathbb{R}_+^2)$ such that $S_{r_1, r_2} \subset \widetilde{O}_{r_1, r_2}$ and $\p_V(S_{r_1, r_2}) = \p_V(\widetilde{O}_{r_1, r_2})$.
Accordingly, \eqref{vxi_2} holds for $\p_V$-almost all $(u_1, u_2)$.

For $n > N(x_1, x_2)$  and $u_i \in [0,\delta]$ one has
$u_i n^{-\beta} < |x_1 - x_2|/2$, $i=1,2$, so
 Lemma \ref{lem_2} applies to $\p(\xi=(r_1,r_2) | V=(u_1,u_2))$. Then
\begin{multline*}
J_{n,k,\beta}(z_1, z_2, u_1, u_2)\! =\! \sum_{r_1, r_2} h(r_1,r_2)
\, \p(\xi_{n,k,1} \!=\! r_1| V\!=\! (u_1,u_2)) \p(\xi_{n,k,2}\! =\! r_2| V\!=\! (u_1,u_2))
 = \mathcal{J}_{n,1} \mathcal{J}_{n,2}
\end{multline*}
where
\vspace{-0.3cm}
\begin{gather*}
\mathcal{J}_{n,i} := \e \left( -\log \left(\frac{\xi_{n,k,i}+1}{k}\right) - H(Y|X) \Big| V=(u_1,u_2)\right)
\\
= \e\left(-\log \left( \frac{\mu_{n,i}+1}{k}\right) - H(Y|X)\right)\p(Y \ne y\big| \|X - x_i\| = u_i n^{-\beta}) \\
+ \e\left(-\log \left( \frac{\mu_{n,i}+2}{k}\right) - H(Y|X)\right)\p\left(Y = y\big| \|X - x_i\| = u_i n^{-\beta}\right) \\
= \e\left(-\log  \left(\frac{\mu_{n,i}+1}{k}\right) + \log f(y_i|x_i)\right)\p\left(Y \ne y\big| \|X - x_i\| = u_i n^{-\beta}\right) \\
+ \e\left(-\log  \left(\frac{\mu_{n,i}+2}{k}\right) + \log f(y_i|x_i)\right)\p\left(Y = y\big| \|X - x_i\| = u_i n^{-\beta}\right) \\
+ (-\log f(y_i|x_i) - H(Y|X))
\end{gather*}
$\mu_{n,i} = \mu_{n}(k, \beta, u_i, x_i, y_i)$, $J_{n,i}= J_{n,i}(k, \beta, u_i, x_i, y_i)$ and $i =1,2$.

According to \eqref{MU} and since $\log k - \log(k-1) \to 0$ as $n \to \infty$
we can write
\begin{gather*}
\sup_{(u,x,y)\in G_n} \left|\e\left(\log  \left(\frac{\mu_{n}(k, \beta, u,x,y)+1}{k}\right)\right) - \log f(y|x)\right| \to 0.
\end{gather*}
Note that $G_n \nearrow (0, \delta] \times \{x \in \mathbb{R}^d \colon f_X(x) > 0 \} \times M$
as $n\to \infty$. For a given version of $f_X$ and any $x\in \mathbb{R}^d$ such that $f_X(x) > 0$, one can find $N(x) \in \mathbb{N}$ to guarantee relation $x \in B_{1, n}$ when $n > N(x)$.
Consider $x_i \in \mathbb{R}^d$ such that
$f_X(x_i)>0$, $i=1,2$. Then, for $n\geq \max\{N(x_1),N(x_2)\}$ and $i=1,2$

\begin{gather*}
\sup_{0<u_i \leq \delta} \left|\log \left( \frac{\mu_{n,i}+1}{k}\right) - \log f(y_i|x_i)\right| \leq \sup_{(u,x,y)\in G_n} \left|\log  \left(\frac{\mu_{n}(k, \beta, u,x,y)+1}{k}\right) - \log f(y|x)\right| \to 0,
\end{gather*}
as $n\to \infty$. In a similar way
\begin{gather*}
\sup_{0<u_i \leq \delta} \left|\log  \left(\frac{\mu_{n,i}+2}{k}\right) - \log f(y_i|x_i)\right| \to 0,\;\;n\to \infty,\;\;i=1,2.
\end{gather*}
Therefore, for $\p_X$-almost all $x_i \in \mathbb{R}^d$ and any $y_i \in M$
($i=1,2$),
\begin{equation}\label{U}
\sup_{0<u_i \leq \delta} |\mathcal{J}_{n,i} - (-\log f(y_i|x_i) - H(Y|X))| \to 0,\;\;n\to \infty.
\end{equation}
Set $F_i := -\log f(y_i|x_i) - H(Y|X)$. Then one has
\begin{gather*}
\sup_{0<u_1, u_2 \leq \delta} |\mathcal{J}_{n,1}\mathcal{J}_{n,2} - F_1 F_2| \leq \sup_{0<u_1 \leq \delta} |\mathcal{J}_{n,1}| \sup_{0<u_2 \leq \delta} |\mathcal{J}_{n,2} - F_2| + |F_2| \sup_{0<u_1 \leq \delta} |\mathcal{J}_{n,1} - F_1|.
\end{gather*}
In view of \eqref{U}, for $n \geq N(x_i, y_i, \delta)$,
the following inequality holds
$\sup_{0<u_1 \leq \delta} |\mathcal{J}_{n,1}| \leq 2|F_1|$.
Whence, for all $y_1, y_2 \in M$ and $\p_X\otimes \p_X$-almost all $(x_1, x_2) \in \mathbb{R}^d\times \mathbb{R}^d$,
$$
\sup_{u_1, u_2 < \delta} |\mathcal{J}_{n,1}\mathcal{J}_{n,2} - F_1 F_2| \to 0,\;\;n\to \infty.
$$
Thus \eqref{ineq7} is proved. It gives us
\begin{multline*}
|\mathcal{U}_{n,k,1}(z_1,z_2) - F_1F_2|
\leq \Big| \int_0^\delta \int_0^\delta J_{n,k,\beta}(z_1, z_2, u_1, u_2) f_{n,k,\beta} (x_1, x_1, u_1, u_2) \, du_1 du_2
- F_1F_2\Big|\\
+ \Big| \int_0^\delta \int_0^\delta J_{n,k,\beta}(z_1, z_2, u_1, u_2) f_{n,k,\beta} (x_1, x_1, u_1, u_2) \, du_1 du_2\Big|  (1 - \p(\rho_{n,k,1}n^{\beta} \leq \delta, \rho_{n,k,2}n^{\beta} \leq \delta)) \\
\leq \sup_{0<u_1, u_2 \leq \delta} |\mathcal{J}_{n,1}\mathcal{J}_{n,2} - F_1 F_2| \int_{0}^{\delta} \int_{0}^{\delta} f_{n,k,\beta}(x_1,x_2,u_1,u_2) \, du_1 du_2 \\
+ |\log k + H(Y|X)| (1 - \p(\rho_{n,k,1}n^{\beta} \leq \delta, \rho_{n,k,2}n^{\beta} \leq \delta)).
\end{multline*}
Moreover,
\begin{gather*}
1 - \p(\rho_{n,k,1}n^{\beta} \leq \delta, \rho_{n,k,2}n^{\beta} \leq \delta) \leq  \p\left( \left\{ \rho_{n,k,1} > \delta n^{-\beta} \right\} \cup \left\{ \rho_{n,k,2} > \delta n^{-\beta} \right\} \right)
\leq 2 \p\left( \rho_{n,k,1} > \delta n^{-\beta} \right).
\end{gather*}
In view of \eqref{f1} and \eqref{ineq4} we get that, for any $z_j\in \mathbb{R}^d\times M$,  $j=1,2$,
$$
|\log k + H(Y|X)| (1 - \p(\rho_{n,k,1}n^{\beta} \leq \delta, \rho_{n,k,2}n^{\beta} \leq \delta)) \to 0,\;\;n\to \infty.
$$
Therefore, $\mathcal{U}_{n,k,1}(z_1,z_2) - F_1F_2\to 0$ as $n\to \infty$.
Hence \eqref{un_ineq} is established and the proof of Theorem~\ref{th_2} is complete. $\square$

{\it Proof of  Corollary 1}.
Let $X \sim N(a, \Sigma)$ where $a \in \mathbb{R}^d$ and $\Sigma > 0$.
It is easily seen that, for any $\varepsilon > 0$, one has $\e |\log f_X(X)|^{\varepsilon} < \infty$ (see, e.g., \cite{BulDim}).
Since, for $x \in \mathbb{R}^d$ and $y \in M$,
\begin{gather*}
f(x, y) = \p(Y = y| X = x) f_X(x), \\
\p(Y = 1| X = x) = \frac{1}{1 + \exp\{-(w, x) - b\}} > 0, \\
f_X(x) = \frac{1}{(2\pi)^{d/2} |\det \Sigma|^{1/2}} \exp\left\{ -\frac{1}{2}(x - a)^T \Sigma^{-1} (x - a) \right\} > 0,
\end{gather*}
we see that $f(\cdot, y)$ is $\mu$-almost everywhere positive.
We show that $f(\cdot, y)$ is $C_0$-constricting for any $y \in \{0, 1\}$. According to Remark 1, it is sufficient to verify that this function is a Lipschitz one. Write
\begin{gather*}
|f(u, y) - f(v, y)| = |\p(Y = y| X = u) f_X(u) - \p(Y = y| X = v) f_X(v)| \\
\leq  \p(Y = y| X = u) |f_X(u) - f_X(v)| + f_X(v) |\p(Y = y| X = u) - \p(Y = y| X = v)| \\
\leq |f_X(u) - f_X(v)| + \max_x |f_X(x)| \, |\p(Y = y| X = u) - \p(Y = y| X = v)|.
\end{gather*}
Note that $\max_x |f_X(x)| < \infty$ and $f_X(\cdot)$ satisfies the Lipschitz condition. Thus it is enough to prove that the function $\p(Y = y|X = x)$ (as a function in $x$) satisfies the Lipschitz condition with a constant $C$.
For any $x \in \mathbb{R}^d$ and $j\in \{1,\ldots,d\}$,
\begin{gather*}
\left|\frac{\partial}{\partial x_j} \p(Y = 1|X = x) \right| = \frac{\exp\{-(w, x) - b\}}{(1 + \exp\{-(w, x) - b\})^2} |w_j| \leq \frac{1}{4} \|w\| < \infty.
\end{gather*}
Thus for $\p(Y = 1|X = x)$ the desired property holds. Obviously, one has $\p(Y = 0|X = x) = 1 - \p(Y = 1|X = x)$ and consequently $\p(Y = 0|X = x)$ satisfies the Lipschitz condition as well.
$\square$

\vskip1cm
{\large {\bf Appendix}}

\vskip0.3cm

{\it Proof of Lemma \ref{lem_1}}.
Fix $n>1$ and
	$k\in \{1,\ldots,n-1\}$.
For $t>0$, $x\in \mathbb{R}^d$ and $\delta \in (0,t)$, introduce an event
$$
	A_{n,k}(x,t,\delta):=\{t-\delta < \rho_{n,k,1}(x) \leq t+\delta\}.
	$$
The event	 $\{\rho_{n,k,1}(x) > t-\delta\}$ means that
there are less than $k$ points (i.e. $0,1,\ldots,k-1$) among $X_2,\ldots,X_n$
in the  ball $B(x,t-\delta)$. The event  $\{\rho_{n,k}(x) \leq t+\delta\}$
signifies that at least $k$ points (i.e. $k,k+1,\ldots,n-1$) among $X_2,\ldots,X_n$
are contained in the  ball $B(x,t+\delta)$. For $x\in \mathbb{R}^d$ and $t>0$, consider
the set
	$$
	S_x(t,\delta):= \{z\in \mathbb{R}^d: t-\delta <\|z-x\|\leq t+\delta\}.
	$$
Let $
	P_x(t,\delta):=\p(X\in S_x(t,\delta)).$
Note that $P_x(t,\delta)\to 0$ as $\delta \to 0+$
since $\mu(S_x(t,\delta))\to 0$ as $\delta \to 0+$, and because $\p_X$ is absolutely continuous w.r.t. the Lebesgue measure  $\mu$.
For $2\leq i_1<\ldots <i_q \leq n$, taking into account the independence of $X_1,\ldots,X_n$, one has
	$$
	\p(X_{i_1}\in S_x(t,\delta),\ldots,X_{i_q}\in S_x(t,\delta))= P_x(t,\delta)^q.
	$$
Note that
	$$
	A_{n,k}(x,t,\delta)= \bigcup_{(s,m)\in J} B_s D_m G_{n-1-(s+m)},
	$$
here $J:=\{(s,m):s\in \{0,\ldots,k-1\}, m\in \{1,\ldots,n-1\}, k\leq s+m\leq n-1\}$, $J=J(n,k)$,
\begin{gather*}
B_s:=\{s \;\mbox{observations among}\;X_2,\ldots,X_n\;\mbox{are in}\;\;B(x,t-\delta)\},\\
D_m:=\{m \;\mbox{observations among}\;X_2,\ldots,X_n\;\mbox{belong to}\;S(t,\delta)\},\\
G_{n-1-(s+m)}:= \{n-1-(s+m) \;\mbox{observations among}\;X_2,\ldots,X_n\;\mbox{are in}\;\mathbb{R}^d\setminus B(x,t+\delta)\}.
\end{gather*}
Clearly, the events $B_s$, $D_m$ and $G_{n-1-(s+m)}$ depend on $x,t,\delta,n$.
We get
\begin{equation}\label{A}
	\p(A_{n,k}(x,t,\delta))= \p(B_{k-1}D_1 G_{n-1-k}) +O(P_x(t,\delta)^2),\;\;\delta\to 0+,
\end{equation}
as $\p(D_m)\!=\! \binom{n-1}{m} P_x(t,\delta)^m (1\! -\! P_x(t,\delta))^{n-1-m}$
and $\p(B_s D_m G_{n-1-(s+m)})\!\leq\! \p(D_m)$ for $(m,s)\in J$.
Set $
	p_x(u)=\p(X\in B(x,u))$, $x\in \mathbb{R}^d$, $u>0$.
Then
	$$
	\p(B_{k-1}D_1 G_{n-1-k})= \binom{n-1}{k-1}(p_x(t-\delta))^{k-1}\binom{n-k}{1} P_x(t,\delta)(1-p_x(t+\delta))^{n-1-k}.
	$$
Indeed, there exist $\binom{n-1}{k-1}$ variants to choose $k-1$ points among $X_2,\ldots,X_n$,
contained in $B(x,t-\delta)$,
and after that $n-1-(k-1)$ variants to choose one point contained in $S_x(t,\delta)$. Other points (their cardinality is $n-1-k$) belong to the complement to the ball $B(x,t+\delta)$.

Now note that, for $r=0,1,\ldots,k-1$,
	$$
	\p(\xi_{n,k,1}(x,y)=r, A_{n,k}(x,t,\delta)) = \p(\xi_{n,k,1}(x,y)=r, B_{k-1}D_1 G_{n-1-k})+ O(P_x(t,\delta)^2),\;\;\delta\to 0+,
	$$
and
	$$
	\p(\xi_{n,k,1}(x,y)=r, B_{k-1}D_1 G_{n-1-k})
	$$
	$$
	= \binom{n-1}{k-1}\binom{k-1}{r} \p(Y=y,X\in B(x,t-\delta))^r \p(Y\neq y,X\in B(x,t-\delta))^{k-1-r}
	$$
	$$
	\times \binom{n-k}{1} P(Y\neq y, X\in S_x(t,\delta))(1-p_x(t+\delta))^{n-1-k}
	$$
\begin{equation}\label{lem_for}
	+ \binom{n-1}{k-1}\binom{k-1}{r-1}\p(Y=y,X\in B(x,t-\delta))^{r-1}\p(Y\neq y,X\in B(x,t-\delta))^{k-1-(r-1)}
\end{equation}
	$$
	\times \binom{n-k}{1} \p(Y=y,X\in S_x(t,\delta))(1-p_x(t+\delta))^{n-1-k}.
	$$
We take into account that there are $\binom{n-1}{k-1}$ variants to choose  $k-1$ points among $X_2,\ldots,X_n$ which lay in $B(x,t-\delta)$, and there exist $\binom{n-k}{1}$ variants to choose among other observations a point $X_q$ belonging to $S_x(t, \delta)$.
Further on there exist two possibilities.
\vskip0.2cm
1. If $Y_q \ne y$ then there are $\binom{k-1}{r}$ variants to choose among points, contained in $B(x,t-\delta)$, $r$ points $X_{i_1},\ldots,X_{i_r}$ such that $Y_{i_m}=y$, $m=1,\ldots,r$. For other $k-1-r$ points $X_{j_1},\ldots,X_{j_{k-1-r}}$, belonging to $B(x,t-\delta)$ one has $Y_{j_s}\neq y$, $s=1,\ldots,k-1-r$.
\vskip0.2cm
2. If $Y_q = y$ then there are $\binom{k-1}{r-1}$ variants to choose among points, contained in $B(x,t-\delta)$, $r-1$ points $X_{i_1},\ldots,X_{i_{r-1}}$ such that $Y_{i_m}=y$, $m=1,\ldots,r-1$. For other $k-1-(r-1)$ points $X_{j_1},\ldots,X_{j_{k-1-(r-1)}}$ belonging to $B(x,t-\delta)$ one has $Y_{j_s}\neq y$, $s=1,\ldots,k-1-(r-1)$.
\vskip0.2cm
Other $n-1-k$ points have to be in the complement of the ball $B(x,t+\delta)$. The probability, for each observation $X_m$,  to be in this complement is equal to $1-p_x(t+\delta)$.
		
For $r=k$, we get
	$$
	\p(\xi_{n,k,1}(x,y)=k, B_{k-1}D_1 G_{n-1-k})
	$$
	$$
	= \binom{n-1}{k-1}\binom{n-k}{1}\p(Y=y,X\in B(x,t-\delta))^{k-1}\p(Y=y,X\in S_x(t,\delta))(1-p_x(t+\delta))^{n-1-k}.
	$$
In this case the reasoning is analogous to the previous one.
The difference is the following. Not only for each  $(k-1)$ points $X_{i_1},\ldots,X_{i_{k-1}}$ (among $X_2,\ldots,X_n$, belonging to $B(x,t-\delta)$),
one has $Y_{i_1}=y,\ldots,Y_{i_{k-1}}=y$, but also for $X_q$ contained in $S_x(t,\delta)$
one has $Y_q=y$. The case $r=k$ is comprised by formula \eqref{lem_for} since $\binom{k-1}{k}=0$.

If a random variable $\tau$ is such that $\tau \geq 0$ a.s., $\e\tau<\infty$ and a random vector $\zeta$ takes values
in $\mathbb{R}^s$ then (see, e.g., \cite{Shiryaev}, Ch. II, Section 7.5)
the function $\e (\tau|\zeta=t)$ can be defined in the following way.
Set $G(B):= \e (\tau\ind \{\zeta\in B\})$ where $B\in \mathcal{B}(\mathbb{R}^s)$. Evidently, $G$
is a finite measure which is absolutely continuous w.r.t.  $\p_{\zeta}$.
Therefore  there is a Borel function $\varphi:\mathbb{R}^s\to \mathbb{R}$ such that, for each $B\in \mathcal{B}(\mathbb{R}^s)$,
$$
\e (\tau\ind\{\zeta\in B\})= \int_{B}\varphi(x)\p_{\zeta}(dx).
$$
In other words $\varphi(t)$ is the Radon - Nikodym derivative $\frac{dG}{d\p_{\zeta}}(t)$, $t\in \mathbb{R}^s$.
Thus $\e(\tau|\zeta)=\varphi(\zeta)$. According to Theorem 5.8.8 \cite{Bogachev} (we take into account that $G\ll\p_{\zeta}$) there exists
\begin{equation}\label{eq_imp}
\lim_{\delta\to 0+}\frac{G(B(t,\delta))}{\p_{\zeta}(B(t,\delta))} = \frac{dG}{d\p_{\zeta}}(t),\;\;t\in \mathbb{R}^s.
\end{equation}
More precisely, this limit exists for $\p_{\zeta}$-almost all $t\in \mathbb{R}^s$ and
is the Radon - Nikodym derivative of the measure $G$ w.r.t. the measure $\p_{\zeta}$, that is a (version) of $\e(\tau|\zeta=t)$.
We employ this result for $\tau=\ind\{\xi_{n,k,1}(x,y)\in D\}$ where $D\in \mathcal{B}(\mathbb{R}_+)$,  $\zeta = \rho_{n,k,1}(x)$, $x\in \mathbb{R}^d$, $y\in M$.
Clearly, $\tau$ is an integrable random variable w.r.t. any finite measure.
Formula \eqref{eq_imp} can be rewritten for $\p_{\rho_{n,k,1}(x)}$-almost all $t\in (0,\infty)$ as follows
\begin{equation}\label{eq_imp_a}
\p(\xi_{n,k,1}(x,y)\in D|\rho_{n,k,1}(x)=t)=\lim_{\delta\to 0+}\frac{\p(\xi_{n,k,1}(x,y)\in D, \rho_{n,k,1}(x)\in B(t,\delta))}{\p(\rho_{n,k,1}(x)\in B(t,\delta))}.
\end{equation}
Note that instead of $B(t,\delta)=[t-\delta,t+\delta]$, where $0<\delta<t$, we can take a set $(t-\delta,t+\delta]$ since, for any $n\in \mathbb{N}$, $n>1$, $k\in \{0,\ldots,n-1\}$ and $x\in \mathbb{R}^d$, a random variable $\rho_{n,k,1}(x)$ has a density.
Indeed, $\p(\rho_{n,k,1}(x) \leq 0)=0$ as there exists a density $f_X(\cdot)$ and, for $t>0$,
\begin{equation}\label{DF}
\p(\rho_{n,k,1}(x) \leq t) = \sum_{j=k}^{n-1}\binom{n-1}{j}p_x(t)^j(1-p_x(t))^{n-1-j}
\end{equation}
where $p_x(t)= \p(X\in B(x,t))$. Evidently,
$p_x(t)=\p(\|X-x\|\leq t)$ is a distribution function of a positive random variable $\|X-x\|$.
At first we show that $p_x(\cdot)$ has a density $f_x(\cdot)$ w.r.t. the Lebesgue measure on $\mathcal{B}(\mathbb{R}_+)$.
After that we prove that there exists a density of a random variable $\rho_{n,k,1}(x)$.

We know that $X=(X_1,\ldots,X_d)$ has a density $f_X(\cdot)$ w.r.t. the Lebesgue measure $\mu$ in $\mathbb{R}^d$ (i.e. $\p_X \ll \mu$). On the other hand, since
$f_X(x)$ is strictly positive for $\mu$-almost all $x\in \mathbb{R}^d$, it is easily seen that $\mu \ll \p_X$.
Consequently, $\p_X\sim \mu$.

Let $\mu_1$ and $\mu_2$ be some measures on a space $(S,\mathcal{B})$ and $h:S\to T$ be
$\mathcal{B}|\mathcal{D}$-measurable function, where $T$ is endowed with a $\sigma$-algebra $\mathcal{D}$.
Introduce the measures $\nu_i:=\mu_i h^{-1}$, $i=1,2$. Then, obviously,
$\mu_1 \ll\mu_2$ yields $\nu_1\ll\nu_2$.
If $Q$ is a Gaussian measure on $\mathcal{B}(\mathbb{R}^d)$
having a density w.r.t. $\mu$, then $Q\sim \mu$ as there exists a strictly
positive version of $dQ/d\mu$ on $\mathbb{R}^d$.
Consider $(S,\mathcal{B}):=(\mathbb{R}^d,\mathcal{B}(\mathbb{R}^d))$,
$(T,\mathcal{D}):=(\mathbb{R}_+,\mathcal{B}(\mathbb{R}_+))$, $h:\mathbb{R}^d\to \mathbb{R}_+$,
where $h(x)=x_1^2+\ldots +x_d^2$ for $x=(x_1,\ldots,x_d)\in \mathbb{R}^d$.
Let $\mu_1=\p_X$ and
$\mu_2$ be a Gaussian law $N(0,I)$ in $\mathbb{R}^d$
with zero mean vector and the unit covariance matrix $I$. Then $\mu_1 \sim \mu_2$
since $\mu_1\sim \mu$ and $\mu_2\sim \mu$. Consequently, $\nu_1 \sim \nu_2$.
Clearly, $\nu_2=\mu_2 h^{-1}$ has the $\chi^2_d$-distribution with a density
w.r.t. the Lebesgue measure $mes_{\mathbb{R}_+}$ on $(\mathbb{R}_+,\mathcal{B}(\mathbb{R}_+))$, i.e.
$$
\p_{\chi^2_d}(u)=\frac{u^{\frac{d}{2}-1}e^{-\frac{u}{2}}}{2^{\frac{d}{2}}
\Gamma(\frac{d}{2})},\;\;u\geq 0.
$$
This density is strictly positive on $(0,\infty)$ and therefore
$\nu_2 \sim mes_{\mathbb{R}_+}$.
Thus $\p_X h^{-1}\sim \nu_2$, hence $\p_X h^{-1}\sim mes_{\mathbb{R}_+}$.
We proved that there exists the density $g$ of a random variable
$X_1^2+\ldots +X_d^2$ w.r.t. the Lebesgue measure on $\mathbb{R}_+$.
Write $d(\p_X h^{-1})/d \,mes_{\mathbb{R}_+}= g$. For any $B\in \mathcal{B}(\mathbb{R}_+)$, one has
$$
\int_B g(t) dt
= \p_Xh^{-1}(B).
$$
If $\int_B g(t) dt=0$ then $\p_X h^{-1}(B)=0$ and hence $mes_{\mathbb{R}_+}(B)=0$.
Take $B:=\{t\in \mathbb{R}_+:g(t)=0\}$. Then $\int_B g(t) dt=0$ and, therefore,
$mes_{\mathbb{R}_+}\{t:g(t)=0\}=0$. In other words, $g$ is strictly positive $mes_{\mathbb{R}_+}$-almost everywhere.

If a random vector  $V$ has a density (w.r.t. measure $\mu$) $q(z)$, $z\in \mathbb{R}^d$,
then, for $x\in \mathbb{R}^d$, the vector $V-x$ has a density $q_x(z)=q(z+x)$, $z\in \mathbb{R}^d$. Consequently, we can
claim that, for each $x\in \mathbb{R}^d$, there exists a density of  random variable $\|X-x\|^2$
w.r.t. the Lebesgue measure $mes_{\mathbb{R}_+}$ on $\mathbb{R}_+$.
This density is strictly positive w.r.t. $mes_{\mathbb{R}_+}$
whenever $f_X(\cdot)$ is strictly positive w.r.t. $\mu$.
If a random variable $\xi_x\geq 0$ has a density  $\gamma_x(u)$, $u\geq 0$ ($x\in \mathbb{R}^d$),
then the random variable $\sqrt{\xi_x}$ has a density
$p_x(u)= 2u \gamma_x(u^2)$, $u\geq 0$.
Thus there is a density $f_x(u)$, $u\geq 0$, of a random variable $\|X-x\|$,
this density is strictly positive for
$mes_{\mathbb{R}_+}$-almost all $u\geq 0$ and $\p_{\|X-x\|}\sim mes_{\mathbb{R}_+}$.

Now we can prove that the density (w.r.t. $mes_{\mathbb{R}_+}$) of a random variable $\rho_{n,k,1}(x)$ has the form
\begin{equation}\label{dens}
h_{n-1,k,x}(u)\!\!=\!\! \sum_{j=k}^{n-1}\!\binom{n-1}{j} (jp_x(u)^{j-1}(1\!-\!p_x(u))^{n\!-\!1\!-\!j}-p_x(u)^j(n\!-\!j\!-\!1)(1\!-\!p_x(u))^{n\!-\!j\!-\!2})f_x(u)
\end{equation}
where $f_x(\cdot)$ is a density corresponding to the distribution function $p_x(\cdot)$, $x\in \mathbb{R}^d$.

It is worth to emphasize that we cannot differentiate the distribution
function to find the density (as the celebrated Cantor function shows).
Thus we have to employ the integral relations. Let $F$ be a distribution function of a positive random variable with a density $f$ (thus $F(0)=0$ and $f(u)=0$, $u<0$). Then using the integration by parts (see, e.g., \cite{Shiryaev}, Ch. II, Section 6.12) and an induction
one can prove that, for each $n\in \mathbb{N}$, a distribution function $F^n(u)$  has a density $nF^{n-1}(u)f(u)$, $u\in \mathbb{R}$.
For $m\in \mathbb{N}$ and $j=0,\ldots,m$, we can write
$$
F(t)^j(1-F(t))^{m-j}
%= \sum_{r=0}^{m-j} F(t)^j C_{m-j}^r 1^r (-F(t))^{m-j-r}
= \sum_{r=0}^{m-j} \binom{m-j}{r} (-1)^{m-j-r} F(t)^{j+m-j-r}
$$
\vspace{-0.3cm}
$$
= \sum_{r=0}^{m-j} \binom{m-j}{r} (-1)^{m-j-r}\int_{(0,t]}(m-r)F^{m-r-1}(u)f(u)du.
$$
The latter formula and \eqref{DF} lead, for $t\geq 0$, to the relation
\begin{equation}\label{z1}
\p(\rho_{n,k,1}(x)\leq t)=\int_0^t h_{n-1,k,x}(u)du
\end{equation}
where $h_{n-1,k,1}$ is given in \eqref{dens}.
We set $h_{n-1,k,1}(u)=0$ for $u<0$.
Note that we come to \eqref{dens} using the polynome of the degree $n-1$ in a distribution
function $p_x(\cdot)$. Hence $h_{n-1,k,x}$ is an integrable function (w.r.t. the Lebesgue measure $mes$ on $\mathbb{R}$).
However, the mentioned polynome has positive and negative coefficients. Therefore, we
have to clarify why $h_{n-1,k,x}$ is a probability density. We explain that if,
for a distribution function $F$, one has
$$
F(t)=\int_{(-\infty,t]}f(u)du,\;\;t\in \mathbb{R},
$$
where $f\in L^1(\mathbb{R},\mathcal{B}(\mathbb{R}),mes)$, then $f$ is a probability density. Clearly, the Lebesgue theorem on dominated convergence yields that $\int_{\mathbb{R}}f(u)du =1$ as $\lim_{t\to \infty}F(t)=1$.
It remains to show that $f(u)\geq 0$ for $mes$-almost all $u$.
Introduce a function
$$
Q(B):= \int_{B}f(u)du,\;\;B\in \mathcal{B}(\mathbb{R}).
$$
Obviously, $Q((a,b])=F(b)-F(a)\geq 0$ for $-\infty \leq a\leq b\leq \infty$
(we set $F(-\infty):=0$, $F(\infty):=1$ and $(a,\infty]:=(a,\infty)$).
Let $G$ be a probability measure on $\mathcal{B}(\mathbb{R})$
generated by a distribution function $F$, then $G((a,b]):=F(b)-F(a)$.
We see that
$G$ and $Q$ coincide on an algebra $\mathcal{A}$ consisting of the finite unions
of the pair-wise disjoint intervals having the form
$(a,b]$, $-\infty \leq a\leq b\leq \infty$. Hence, $Q$ is a finite
nonnegative function on $\mathcal{A}$.
Clearly, $Q$ is a countably additive function on $\mathcal{B}(\mathbb{R})$
and $Q(\mathbb{R})$ is finite.
It remains to note that, for any $B\in \mathcal{B}(\mathbb{R})$
and each $\varepsilon>0$, there exists $A\in \mathcal{A}$ such that $|Q(B)-Q(A)|<\varepsilon$.
Indeed,
$$
|Q(B)-Q(A)|=\left|\int_{\mathbb{R}}(\ind\{B\}-\ind\{A\})f(u)du\right|
\leq
\int_{\mathbb{R}}\ind\{B\triangle A\}|f(u)|du.
$$
Consequently, for each $B\in \mathcal{B}(\mathbb{R})$, one can find $A_n\in \mathcal{A}$ ($n\in \mathbb{N}$) such that $Q(A_n)\to Q(B)$ as $n\to \infty$. Taking into account that $Q(A_n)\geq 0$ we get $Q(B)\geq 0$.
Assume now that $\mu(E)>0$ where $E=\{x:f(x)<0\}$.
Note that $E=\cup_{m=1}^{\infty}\{-\infty<f(x)\leq -\frac{1}{m}\}$.
Then in a standard way we come to the contradiction. Therefore, $\mu(E)=0$.
Thus formula \eqref{dens} provides a probability density
$h_{n-1,k,x}(\cdot)$ of
the random variable $\rho_{n,k,1}(x)$ distribution where $x\in \mathbb{R}^d$.

Hence, for each $x\in \mathbb{R}^d$, $y\in M$, $r\in \{0,1,\ldots,k\}$ and $\p_{\rho_{n,k,1}(x)}$-almost all $t\in (0,\infty)$ in view of \eqref{eq_imp_a} one has
\begin{equation}\label{c_rho}
	\p(\xi_{n,k,1}(x,y)=r|\rho_{n,k,1}(x)=t)
	= \lim_{\delta\to 0+}\frac {\p(\xi_{n,k,1}(x,y)=r,A_{n,k}(x,t,\delta))}{\p(A_{n,k}(x,t,\delta))}.
\end{equation}
Applying the expressions obtained for nominator and denominator of the latter fraction
in \eqref{c_rho}
and taking into account that a function
$\p(X\in B(x,t))$ is continuous in $(x,t)\in \mathbb{R}^d\times \mathbb{R}_+$ (see, e.g., Lemma 1 in \cite{BulDim}) we get, for each $n\in \mathbb{N}$ ($n>1$), $k \in \{1,\ldots,n-1\}$, $r=0,\ldots,k$ and $\p_{\rho_{n,k,1}(x)}$-almost all $t\in (0,\infty)$,
	$$
	\p(\xi_{n,k,1}(x,y)=r|\rho_{n,k,1}(x)=t)
	$$
\begin{equation}\label{eq_nor}
	= \binom{k-1}{r} p^r(1-p)^{k-1-r} \lim_{\delta\to 0+}\frac{\p(Y\neq y,X\in S_x(t,\delta))}{\p(X\in S_x(t,\delta))}
\end{equation}
$$
+
	\binom{k-1}{r-1}p^{r-1}(1-p)^{k-r}\lim_{\delta\to 0+}\frac{\p(Y= y,X\in S_x(t,\delta))}{\p(X\in S_x(t,\delta))}
$$
where $p:=\p(Y=y|X\in B(x,t))$, $p=p(x,y,t)$. We used that $\p(X\in B(x,t))>0$, for
$\mu$-almost all $x\in \mathbb{R}^d$ and $t>0$, since $f_X(z)$
is strictly positive for $\mu$-almost all $z\in \mathbb{R}^d$.
However, we have to explain the existence of limits in \eqref{eq_nor}.
Let us employ formula \eqref{eq_imp}
for $\zeta:= \|X-x\|$, $G(D):= \e (\ind\{\zeta \in D\}\ind\{Y=y\})$, $x\in \mathbb{R}^d$, $y\in M$ and $D\in \mathcal{B}(\mathbb{R})$. We can claim that, for $x\in \mathbb{R}^d$ and $y\in M$,
the limits appearing in \eqref{eq_nor} exist for $\p_{\|X-x\|}$-almost all $t>0$.
Indeed,
\begin{equation}\label{alpha}
	\lim_{\delta\to 0+}\frac{\p(Y\neq y,X\in S_x(t,\delta))}{\p(X\in S_x(t,\delta))}
= \lim_{\delta\to 0+}\frac{\p(Y\neq y,t-\delta<\|X-x\|\leq t+\delta)}{\p(t-\delta<\|X-x\|\leq t+\delta)}
\end{equation}
$$
	=\p(Y\neq y|\|X-x\|=t):=\alpha(x,y,t).
	$$
We have seen that $\p_{\|X-x\|}\sim mes_{\mathbb{R}_+}$ for each $x\in \mathbb{R}^d$.
Therefore
$\p(X\in S_x(t,\delta))>0$ for all $x\in \mathbb{R}^d$ and $\delta>0$.
Moreover, for each $y$ belonging to a finite set $M$ and $x\in \mathbb{R}^d$, the limits in
\eqref{eq_nor} exist for $mes_{\mathbb{R}_+}$-almost all $t\in (0,\infty)$ as
$\p_{\|X-x\|}\sim mes_{\mathbb{R}_+}$.
Consequently, the measure $\p_{\rho_{n,k,1}(x)}$ of a set of points $t\in \mathbb{R}_+$ such that the limits in \eqref{eq_nor} do not exist equals to zero since $\p_{\rho_{n,k,1}(x)} \ll mes_{\mathbb{R}_+}$.
The proof is complete. $\square$

\vskip0.2cm

	{\it Proof of Lemma \ref{lem_2}}.	
Fix arbitrary $z_j=(x_j,y_j)\in \mathbb{R}^d \times M$, $j\in \{1,2\}$, such that $x_1\neq x_2$. Consider $n\in \mathbb{N}$, $n>2$, and $k\in \{1,\ldots,n-1\}$. Note that
	\begin{gather*}
	\xi_{n,k,1}(z_1,z_2) = \sharp\{i: i\in \{3,\ldots,n\}, Y_i=y_1,  X_i\in B(x_1,R_1)\} + \ind\{ y_2 = y_1,  x_2 \in B(x_1,R_1)\}, \\
	\xi_{n,k,2}(z_1,z_2) = \sharp\{i: i\in \{3,\ldots,n\}, Y_i=y_2,  X_i\in B(x_2, R_2)\} + \ind\{ y_1 = y_2,  x_1 \in B(x_2,R_2)\},
	\end{gather*}
here  $R_j:=\rho_{n,k,j}(x_1,x_2)$, $j=1,2$. Recall that $\rho_{n,k,j}(x_1,x_2)$ and $\xi_{n,k,j}(z_1,z_2)$ are
defined in \eqref{ro} and \eqref{hi}, respectively.
Take any $\varepsilon \in (0,|x_1 - x_2|/2)$ and $t_j\in (0,|x_1 - x_2|/2 - \varepsilon) $, $j=1,2$. Then there exists $\delta >0$ such that $\delta < t_j$,  $j=1,2$, and $t_1 + t_2 + 2\delta < |x_1 - x_2|$. Introduce the events
$$
	A_{n,k} :=\bigcap_{j=1}^2 \{t_j-\delta < \rho_{n,k,j}(x_1, x_2) \leq t_j+\delta\}.
$$
Clearly,   $A_{n,k} = A_{n,k}(x_1, x_2, t_1, t_2, \delta)$.
Further on in the proof we consider $j \in \{1, 2\}$ without mentioning.
To simplify the exposition we use a notation similar in meaning to that employed for proving Lemma \ref{lem_1}. However, we have to emphasize that now we use vectors with two components in contrast to random variables appearing in the proof of Lemma \ref{lem_1}.
For instance, $A_{n,k}$  is not the same as previously.
An event $\{\rho_{n,k,j}(x_1, x_2) > t_j-\delta\}$ means that the closed ball $B(x_j, t_j-\delta)$ contains less than $k$ (i.e. $0,1,\ldots,k-1$)
observations among $X_3, \ldots,X_n$, since $x_i$, $i\in \{1,2\}\setminus \{j\}$, does not belong to this ball as
$|x_1 - x_2| > t_j - \delta$.
An event $\{\rho_{n,k,j}(x_1, x_2) \leq t_i+\delta\}$ signifies that in $B(x_j,t_j+\delta)$ there are at least
	$k$ (i.e. $k,k+1,\ldots,n-2$) points among  $X_3,\ldots,X_n$ because $x_i$, $i\in \{1,2\}\setminus \{j\}$,
does not belong to this ball as $|x_1 - x_2| > t_j + \delta$.
One has
$$
	A_{n,k} = \bigcup_{(s_1,s_2,m_1,m_2)\in J_{n,k}}
	B_{1,s_1}D_{1,m_1}B_{2,s_2}D_{2,m_2} G_{n-2 - (s_1+s_2+m_1+m_2)}
$$
where $J_{n,k}$ consists of $(s_1,s_2,m_1,m_2)$ such that
$s_1, s_2 \!\in\! \{0,\ldots,k-1\}$, $m_1, m_2 \!\in\! \{1,\ldots,n-1\}$, $s_1+m_1\geq k$, $s_2+m_2\geq k$, $s_1+s_2+m_1+m_2 \leq n-2$ and
	$$
	B_{j,s}:=\{s \;\mbox{variables among}\;X_3,\ldots,X_n\;\mbox{belong to}\;B(x_j,t_j-\delta)\},
	$$
	$$
	D_{j,m}:=\{m \;\mbox{variables among}\;X_3,\ldots,X_n\;\mbox{belong to}\;S_{x_j}(t_j,\delta)\},
	$$
	$$
	G_{l}:= \left\{l \;\mbox{variables among}\;X_3,\ldots,X_n\;\mbox{belong to}\;\mathbb{R}^d\setminus \bigcup_{j=1}^2 B(x_j, t_j+\delta)\right\}.	
	$$
Since $k \leq [(n-2)/2]$ the set $J_{n,k}$ is nonempty as $(k-1, k-1, 1, 1) \in J_{n,k}$.
More precisely, $B_{j,s}=B_{j,s}(x_j,t_j,\delta,n)$, $D_{j,m}=D_{j,m}(x_j,t_j,\delta,n)$ and
$G_{l} = G_{l}(x_1, x_2, t_1, t_2, \delta)$.
In a way similar to \eqref{A} one has
	$$
	\p(A_{n,k})= \p(B_{1,k-1}D_{1,1}B_{2,k-1}D_{2,1} G_{n-2-2k}) +O(P_{x_1}(t_1,\delta)^2P_{x_2}(t_2,\delta)) + O(P_{x_1}(t_1,\delta)P_{x_2}(t_2,\delta)^2)
	$$
as $\delta\to 0+$ because, for $(s_1,s_2,m_1,m_2)\in J_{n,k}$,
	$$
	\p(D_{1,m_1} D_{2,m_2} )\!=\! \binom{n-2}{m_1}  P_{x_1}(t_1,\delta)^{m_1} \binom{n\!-\!2\!-\!m_1}{m_2}  P_{x_2}(t_2,\delta)^{m_2} (1 \!-\! P_{x_1}(t_1,\delta) - P_{x_2}(t_2,\delta))^{n\!-\!2\!-\!m_1\!-\!m_2},$$
	$$
	\p(B_{1,s_1}D_{1,m_1}B_{2,s_2}D_{2,m_2} G_{n-2 - (s_1+s_2+m_1+m_2)})\leq \p(D_{1,m_1} D_{2,m_2}),
	$$
$P_x(t,\delta)=\p(X\in S_x(t,\delta))$, $x\in \mathbb{R}^d$, $t>0$ and $\delta >0$. It is easily seen that
	\begin{gather*}
	\p(B_{1,k-1}D_{1,1}B_{2,k-1}D_{2,1} G_{n-2-2k})
	= Poly(k-1, 1, k-1, 1, n-2-2k) \\
	\times (p_{x_1}(t_1-\delta))^{k-1} P_{x_1}(t_1,\delta)(p_{x_2}(t_2-\delta))^{k-1} P_{x_2}(t_2,\delta) (1-p_{x_1}(t_1+\delta) - p_{x_2}(t_2+\delta))^{n-2-2k}
	\end{gather*}
where $p_x(t)=\p(X\in B(x,t))$, $x\in \mathbb{R}^d$, $t>0$ and
$$
Poly(k_1,\ldots,k_q):= \frac{(k_1+\ldots +k_q)!}{k_1!\ldots k_q!}, \;\;k_i\in \mathbb{Z}_+,\;i=1,\ldots,q.
$$	
Thus, for $r_1, r_2 \in \{0, \ldots, k\}$,
	\begin{gather*}
	\p(\xi_{n,k,1}(z_1,z_2) = r_1, \xi_{n,k,2}(z_1,z_2) = r_2, A_{n,k}) \\
	= \p(\xi_{n,k,1}(z_1,z_2) = r_1, \xi_{n,k,2}(z_1,z_2) = r_2, B_{1,k-1}D_{1,1}B_{2,k-1}D_{2,1} G_{n-2-2k})\\+ O(P_{x_1}(t_1,\delta)^2P_{x_2}(t_2,\delta)) + O(P_{x_1}(t_1,\delta)P_{x_2}(t_2,\delta)^2),\;\;\delta\to 0+.
	\end{gather*}
Introduce the auxiliary events.
Let $B_{j,s}^l$ mean that $s$ observations among $X_3,\ldots,X_n$ are contained in $B(x_j,t_j-\delta)$ while the rest are not, moreover,
$l$ points among $X_i$'s contained in this ball, i.e. $X_{i_1},\ldots,X_{i_l}$ are such that $Y_{i_m}=y_j$, $m=1,\ldots,l$. Clearly, $B_{j,s}^l=B_{j,s}^l(x_j,y_j,t_j,\delta,n)$.
Analogously one can define an event $D_{j,s}^l$ (namely, $s$ points among $X_3,\ldots,X_n$ are in
$S_{x_j}(t_j,\delta)$ and other ones do not belong to this set, moreover, $l$ points among $X_i$'s belonging to $S_{x_j}(t_j,\delta)$ are such that corresponding $Y_i=y_j$). Note that
$D_{j,m}^l=D_{j,m}^l(x_j,y_j,t_j,\delta,n)$.
Then, for $r_1, r_2 \in \{0,1,\ldots,k\}$,
$$
	\{\xi_{n,k,1}(z_1,z_2) = r_1, \xi_{n,k,2}(z_1,z_2) = r_2\} \cap B_{1,k-1}D_{1,1}B_{2,k-1}D_{2,1} G_{n-2-2k} $$
$$
	= B_{1,k-1}^{r_1} D_{1,1}^0 B_{2,k-1}^{r_2} D_{2,1}^0 G_{n-2-2k} \cup B_{1,k-1}^{r_1-1} D_{1,1}^1 B_{2,k-1}^{r_2} D_{2,1}^0 G_{n-2-2k}
$$
\begin{equation}\label{uni}
\cup	B_{1,k-1}^{r_1} D_{1,1}^0 B_{2,k-1}^{r_2-1} D_{2,1}^1 G_{n-2-2k} \cup
	B_{1,k-1}^{r_1-1} D_{1,1}^1 B_{2,k-1}^{r_2-1} D_{2,1}^1 G_{n-2-2k}.
\end{equation}
If $r=0$ then $B_{j,k-1}^{r-1}:=\varnothing$ (for $j\in \{1,2\}$).
Evidently, four events appearing in the union in \eqref{uni} are pair-wise disjoint.
We evaluate their probabilities.	
One has
	\begin{gather*}
	\p(B_{1,k-1}^{r_1} D_{1,1}^0 B_{2,k-1}^{r_2} D_{2,1}^0 G_{n-2-2k})
	= \\
	Poly(k-1, 1, k-1, 1, n-2-2k) \binom{k-1}{r_1} \binom{k-1}{r_2}
	(1-p_{x_1}(t_1+\delta) - p_{x_2}(t_2+\delta))^{n-2-2k}\\
	\times  \p(Y\!=\!y_1,X\!\in\! B(x_1,t_1-\delta))^{r_1}
	\p(Y\!\neq\! y_1,X\!\in\! B(x_1,t_1-\delta))^{k-1-r_1} P(Y\!\neq\! y_1, X\in S_{x_1}(t_1,\delta))\\
	\times \p(Y\!=\!y_2,X\!\in\! B(x_2,t_2-\delta))^{r_2}
	\p(Y\neq y_2, X\!\in\! B(x_2,t_2-\delta))^{k-1-r_2} P(Y\!\neq\! y_2, X\in S_{x_2}(t_2,\delta)).
	\end{gather*}
Indeed, there are $Poly(k-1, 1, k-1, 1, n-2-2k)$ variants for partitioning of $X_3,\ldots,X_n$
into groups belonging, correspondingly, to pair-wise disjoint (under conditions imposed on
$t_1,t_2,|x_1-x_2|$ and $\delta$) sets $B(x_1,t_1-\delta)$, $S_{x_1}(t_1,\delta)$, $B(x_2,t_2-\delta)$, $S_{x_2}(t_2,\delta)$ and
$\mathbb{R}^d\setminus \cup_{j=1}^2 B(x_j,t_j+\delta)$.
We note that there exist 	$\binom{k-1}{r_1}$ variants to choose ${r_1}$
points $X_i$, $i\in I$, among $X_{q_1},\ldots, X_{q_{k-1}}$ ($3\leq q_1<\ldots < q_{k-1}\leq n$) belonging to $B(x_1,t_1-\delta)$ such that  $Y_i=y_1$ for $i\in I$
and $Y_q\neq y_1$ for $q\in \{q_1,\ldots,q_{k-1}\}\setminus I$,
$\sharp I=r_1$.
In a similar way one can explain the appearance of a factor $\binom{k-1}{r_2}$.
For other three events their probabilities can be found analogously. As a result we obtain
	\begin{gather*}
	\p(\xi_{n,1} = r_1, \xi_{n,2} = r_2, B_{1,k-1}D_{1,1}B_{2,k-1}D_{2,1} G_{n-2-2k}) \\
	= Poly(k-1, 1, k-1, 1, n-2-2k)
	(1-p_{x_1}(t_1+\delta) - p_{x_2}(t_2+\delta))^{n-2-2k}\\
	\times \prod_{j=1}^2  \Biggl( \binom{k-1}{r_j}\p(Y=y_j,F(x_j,t_j,\delta))^{r_j}
	\p(Y\neq y_j,F(x_j,t_j,\delta))^{k-1-r_j} P(Y\neq y_j, X\in S_{x_j}(t_j,\delta))\\
	+ \binom{k-1}{r_j-1}\p(Y=y_j,F(x_j,t_j,\delta))^{r_j-1}
	\p(Y\neq y_j,F(x_j,t_j,\delta))^{k-r_j} P(Y = y_j, X\in S_{x_j}(t_j,\delta))\Biggr)
	\end{gather*}
where $F(x_j,t_j,\delta):=\{X\in B(x_j,t_j-\delta)\}$, $j=1,2$.	
Hence, in view of \eqref{alpha} and as, for each $x\in \mathbb{R}^d$, the distribution of $\|X-x\|$ is equivalent on
$\mathbb{R}_+$ to the Lebesgue measure $mes_{\mathbb{R}_+}$, we can state the following. For any $x_1,x_2\in \mathbb{R}^d$ ($x_1\neq x_2$), and $mes\otimes mes$-almost all $t=(t_1,t_2) \in (0,\infty)\times (0,\infty)$
such that $t_1+t_2 <|x_1-x_2|$, one has
\begin{equation}\label{A}
	\lim_{\delta\to 0+}\frac{\p(\xi_{n,k,1}(z_1,z_2) = r_1, \xi_{n,k,2}(z_1,z_2) = r_2, A_{n,k})}{\p(A_{n,k})}
\end{equation}
$$
	=\lim_{\delta\to 0+}\prod_{j=1}^2  \Bigg\{\! \binom{k-1}{r_j}\!
	\left(
	\frac{\p(Y\!=\!y_j,X\!\in\! B(x_j,t_j-\delta))}
	{p_{x_1}(t_1\!-\!\delta)}
	\right)^{r_j}\!
	\left(
	\frac{\p(Y\!\neq\! y_j,X\!\in\! B(x_j,t_j\!-\!\delta))}{p_{x_1}(t_1\!-\!\delta)}
	\right)^{k-1-r_j}
$$
$$
 \times
	\left(
	\frac{P(Y\!\neq\! y_j, X\!\in\! S_{x_j}(t_j,\delta))}
	{P(X\!\in\! S_{x_j}(t_j,\delta))}
	\right)
$$
$$
	+ \binom{k-1}{r_j-1}
	\!\left(
	\frac{\p(Y\!=\!y_j,X\!\in\! B(x_j,t_j\!-\!\delta))}
	{p_{x_1}(t_1\!-\!\delta)}
	\right)^{r_j-1}
	\!\left(
	\frac{\p(Y\!\neq\! y_j,X\!\in\! B(x_j,t_j\!-\!\delta))}{p_{x_1}(t_1\!-\!\delta)}
	\right)^{k-r_j}
$$
$$
\times
	\left(
	\frac{P(Y\! =\! y_j, X\!\in\! S_{x_j}(t_j,\delta))}
	{P(X\!\in\! S_{x_j}(t_j,\delta))}
	\right)
	\Bigg\}
$$
$$
	=\prod_{j=1}^2  \Bigg\{ \binom{k-1}{r_j}
	p_j^{r_j}
	(1 - p_j)^{k-1-r_j}
	\lim_{\delta \to 0+}
	\left(
	\frac{P(Y\neq y_j, X\in S_{x_j}(t_j,\delta))}
	{P(X\in S_{x_j}(t_j,\delta))}
	\right)
$$
$$
	+ \binom{k-1}{r_j-1}
	p_j^{r_j-1}
	(1-p_j)^{k-r_j}
	\lim_{\delta \to 0+}
	\left(
	\frac{P(Y = y_j, X\in S_{x_j}(t_j,\delta))}
	{P(X\in S_{x_j}(t_j,\delta))}
	\right)
	\Bigg\}
$$
$$
=\prod_{j=1}^2  \Bigg\{ \binom{k-1}{r_j}
	p_j^{r_j}
	(1 - p_j)^{k-1-r_j}
	\alpha(x_j,y_j,t_j) + \binom{k-1}{r_j-1}
	p_j^{r_j-1}
	(1-p_j)^{k-r_j}(1-\alpha(x_j,y_j,t_j))\Bigg\}
$$
where $p_j := p_j(x_j, y_j, t_j) = \p(Y=y_j|X\in B(x_j,t_j))$ and we use that
a function $\p(X\in B(x,t))$ is continuous in $(x,t)\in \mathbb{R}^d\times \mathbb{R}_+$.
We employ also that $\p(X\in B(x,t))>0$, for
all $x\in \mathbb{R}^d$ and $t>0$, since $f_X(z)$ is strictly positive for $\mu$-almost all $z\in \mathbb{R}^d$.

Note that the proof of Theorem 5.8.8 in \cite{Bogachev} shows that formula \eqref{eq_imp}
is also valid if we replace in it the balls $B(t,\delta)$ by the cubes $\widetilde{B}(t,\delta):=
\prod_{i=1}^s [t_i-\delta,t_i+\delta]$ where $t=(t_1,\ldots,t_s)\in \mathbb{R}^s$. Thus, for $\p_{\zeta}$-almost all $t\in \mathbb{R}^s$,
one can write instead of \eqref{eq_imp} that
\begin{gather*}
\lim_{\delta\to 0+}\frac{G(\tilde{B}(t,\delta))}{\p_{\zeta}(\tilde{B}(t,\delta))} = \e(\tau|\zeta = t)
\end{gather*}
Take now $s=2$, $
\tau \!:=\! \mathbb{I}\{\xi_{n,k,1}(z_1,z_2)\! =\! r_1, \xi_{n,k,2}(z_1,z_2) = r_2 \},\;\;
\zeta\! :=\! ((\rho_{n,k,1}(x_1,x_2),\rho_{n,k,2}(x_1,x_2))
$ where  $n\in \mathbb{N}$, $n>2$, $k=0,\ldots,n-1$, $z_1=(x_1,y_1)\in \mathbb{R}^d\times M$, $z_2=(x_2,y_2)\in \mathbb{R}^d\times M$ and $r_1,r_2\!\in\! \{0,\ldots,k\}$.
Then, for $\p_{\zeta}$-almost all $t = (t_1, t_2)\in (0,\infty)\times (0,\infty)$, we get
$$
\p(\xi_{n,k,1}(z_1,z_2)=r_1, \xi_{n,k,2}(z_1,z_2)=r_2|\rho_{n,k,1}(x_1, x_2) = t_1, \rho_{n,k,2}(x_1, x_2) = t_2)
$$
\begin{equation}\label{E}
= \lim_{\delta \to 0+} \frac{\e(\mathbb{I}\{ \zeta \in \widetilde{B}(t, \delta) \} \tau)}{\p(\zeta \in \widetilde{B}(t, \delta))}.
\end{equation}

Let us show that, for each $x_1,x_2\in \mathbb{R}^d$ ($x_1\neq x_2$) and $\p_{\zeta}$-almost all $t = (t_1, t_2)$
belonging to the set $B_{x_1, x_2}(\varepsilon) = \{t \in \mathbb{R}^2: 0 < t_j \leq |x_1 - x_2|/2 - \varepsilon, \,\, j =1, 2\}$, the latter limit coincides with the obtained value for $\lim_{\delta\to 0+}\p(\xi_{n,k,1}(z_1,z_2)\!=\!r_1, \xi_{n,k,2}(z_1,z_2)\!=\!r_2|A_{n,k})$ in \eqref{A}
where $A_{n,k}$ depends on $x_1,x_2$ and $\delta$.

For this purpose we demonstrate that, for any $x_1,x_2\in \mathbb{R}^d$ ($x_1\neq x_2$), there exists a
positive
measurable function $f_{n,k}(x_1,x_2,\cdot,\cdot)$ such that if $B\in \mathcal{B}(\mathbb{R}^2)$ and $B\subset B_{x_1, x_2}(\varepsilon)$ then
\begin{equation}\label{D}
\p(\zeta \in B)= \int_B f_{n,k}(x_1,x_2,u_1,u_2)du_1 du_2.
\end{equation}
For $(u_1,u_2)\in B_{x_1, x_2}(\varepsilon)$, one can write
\begin{multline*}
\p(\rho_{n,k,1}(x_1,x_2) \leq u_1, \rho_{n,k,2}(x_1,x_2) \leq u_2) \\
= \sum_{(l_1,l_2)\in J(n,k)} \binom{n}{l_1} \binom{n-l_1}{l_2} p_{x_1}^{l_1}(u_1) p_{x_2}^{l_2}(u_2)
(1 - p_{x_1}(u_1) - p_{x_2}(u_2)))^{n-l_1-l_2}
\end{multline*}
where $J(n,k):=\{(l_1,l_2):l_1 \geq k, l_2 \geq k, l_1+l_2 \leq n-2\}$.
This set is nonempty as  $(k, k) \in J(n, k)$.

Consequently,
we get a polynome in $p_{x_1}(u_1)$ and $p_{x_2}(u_2)$. In the proof of Lemma \ref{lem_1} we have seen that, for $x\in \mathbb{R}^d$, $u>0$ and $l\in \mathbb{N}$, the distribution function $p_x(u)^l$ has a density $l p_x^{l-1}(u)f_x(u)$ where $f_x(\cdot)$ is a density
of $p_x(\cdot)$ w.r.t. $mes_{\mathbb{R}_+}$. Thus $p_{x_1}^{l_1}(u_1) p_{x_2}^{l_2}(u_2)$ is absolutely continuous
w.r.t. $mes_{\mathbb{R}_+}\otimes mes_{\mathbb{R}_+}$. Hence there is an integrable
 (w.r.t. the restriction of $mes \otimes mes$ on $B_{x_1,x_2}(\varepsilon)$) function $f_{n,k}(x_1,x_2,\cdot,\cdot) \colon B_{x_1, x_2}(\varepsilon) \to \mathbb{R}$ such that, for a set $(0, u_1] \times (0, u_2]$, $(u_1,u_2) \in B_{x_1, x_2}(\varepsilon)$, formula \eqref{D} takes place. The additivity of the integral implies the validity of the mentioned formula for any
parallelepiped $(a_1, b_1] \times (a_2, b_2]$ where $a_i \leq b_i$ ($i = 1,2$), $(a_1, a_2), (b_1, b_2) \in B_{x_1, x_2}(\varepsilon)$. Moreover, formula \eqref{D} holds for an algebra $\mathcal{E}$ of subsets of $B_{x_1, x_2}(\varepsilon)$ which can be represented as a finite union
of such parallelepipeds. Thus we have seen that
$$
Q(B) := \int_B f_{n,k}(x_1,x_2,u_1,u_2)\,du_1 du_2,\;\;B\in \mathcal{B}(\mathbb{R}^2)\cap B_{x_1,x_2}(\varepsilon),
$$
and $G(B) := \p(\zeta \in B)$ coincides on $\mathcal{E}$.
In a similar way to the proof of Lemma \ref{lem_1} we get that
$
f_{n,k}(x_1,x_2,u_1,u_2) \geq 0
$
for $(mes \otimes mes)$-almost all $u = (u_1, u_2)\in B_{x_1, x_2}(\varepsilon)$.
Therefore, the desired formula \eqref{D} is established.

Compare \eqref{A} and \eqref{E}.
We show now that, for $(t_1,t_2)\in B_{x_1,x_2}(\varepsilon)$ and all $\delta >0$ small enough (i.e. for $\delta < \Delta(x_1,x_2,t_1,t_2)$), one has $\e(\ind\{\zeta \in C(t,\delta)\}\tau)=\e(\ind\{\zeta \in \widetilde{B}(t,\delta)\tau)$ and
$\p(\zeta \in C(t,\delta))=\p(\zeta \in \widetilde{B}(t,\delta))$ where
$C(t,\delta):=(t_1-\delta,t_1+\delta]\times (t_2-\delta,t_2+\delta]$.
Using the relation
$C(t,\delta)\subset \widetilde{B}(t,\delta)\subset B_{x_1, x_2}(\varepsilon)$ for all $\delta \in (0,\Delta(x_1,x_2,t_1,t_2))$ and due to \eqref{D}, we get
$\p(\zeta \in \widetilde{B}(t,\delta)\setminus C(t,\delta))=0$ as
$(mes \otimes mes) (\widetilde{B}(t,\delta)\setminus C(t,\delta))=0$.
Taking into account that the restriction of $\p_{\zeta}$ to $B_{x_1,x_2}(\varepsilon)$ is absolutely continuous
w.r.t. the corresponding restriction of $mes_{\mathbb{R}_+}\otimes mes_{\mathbb{R}_+}$ we can claim that,
for each $x_1,x_2 \in \mathbb{R}^d$ ($x_1\neq x_2$) and for $\p_{\zeta}$-almost all
$(t_1,t_2)\in B_{x_1,x_2}(\varepsilon)$, formulas \eqref{conind} and \eqref{claw} are established.

Note that $B_{x_1, x_2}(\varepsilon) \uparrow (0, |x_1 - x_2|/2) \times (0, |x_1 - x_2|/2)$ as $\varepsilon \to 0$. Consequently, formulas \eqref{conind} and \eqref{claw} are valid $\p_\zeta$-a.s. for points of the set $(0, |x_1 - x_2|/2) \times (0, |x_1 - x_2|/2)$.
$\square$

{\it Proof of Lemma 3}.
One has
$$
\e(W\ind\{V\in B\})=\e(\e(W\ind\{V\in B\})|V)= \e(\ind\{V\in B\}\e(W|V))
$$
$$
= \int_{\mathbb{R}^m}\ind\{t\in B\}\e(W|V=t)\p_V(dt)= \int_B\e(W|V=t)\p_V(dt).
$$
Consequently,
$$
\e(W|V\in B)= \frac{\e(W\ind\{V\in B\})}{\p(V\in B)}=\frac{\int_B\e(W|V=t)\p_V(dt)}{\p(V\in B)}=\int_B\e(W|V=t)\widetilde{\p}_{V,B}(dt)
$$
where $\widetilde{\p}_{V,B}(D)= \frac{\p(D\cap B)}{\p(B)}$, $D\in \mathcal{B}(\mathbb{R}^m)$. $\square$

{\it Proof of Lemma 4}. It is enough to demonstrate that, for any set $B \in \mathcal{B}(\mathbb{R})$,
$$
\int_B \e(\xi|\eta = t) \, \p_\eta(dt) = \sum_{r \in S} \int_B \e(\xi|\zeta = r, \eta=t) \p(\zeta=r|\eta=t) \, \p_\eta(dt).
$$
Clearly,
$$
\int_B \e(\xi|\eta = t) \, \p_\eta(dt) = \e \xi \mathbb{I}(\eta \in B) = \sum_{r \in S} \e \xi \mathbb{I}(\eta \in B, \zeta = r)
$$
$$
= \sum_{r \in S} \int_{B \times \{r\} } \e(\xi|\eta=t, \zeta = v) d\p_{\eta,\zeta}(t,v)
$$
where $d\p_{\eta,\zeta}(t,r)$ means the integration w.r.t.
measure $\p(\eta\in \cdot,\zeta\in \cdot)$.
Now we show that, for a measurable function $\varphi \colon \mathbb{R} \times S \to \mathbb{R}$ such that
$\e |\varphi(\eta, \zeta)| < \infty$, the following relation holds
\begin{equation}\label{V}
\int_{B \times \{r\} } \varphi(t,v) d\p_{\eta,\zeta}(t,v) =  \int_{B} \varphi(t, r) \p(\zeta = r|\eta = t)\p_{\eta}(dt).
\end{equation}
Indeed, for $A\in \mathcal{B}(\mathbb{R})$ and $s\in S$, consider
$\varphi(t, v) := \ind\{t \in A, v = s\}$, $t\in \mathbb{R}$ and $v\in S$. Obviously, if $ s \ne r$ then \eqref{V} is true. If $s = r$ then
\begin{gather*}
\int_{B \times \{r\} } \varphi(t,v) d\p_{\eta, \zeta}(t,v) = \p_{\eta,\zeta}((A \cap B) \times \{r\}) = \p(\zeta = r, \eta \in A \cap B)
\\
=
\int_{A \cap B} \p(\zeta = r| \eta = t) \p_{\eta}(dt) = \int_{B} \mathbb{I}(t \in A, r = r) \p(\zeta = r| \eta = t) \p_{\eta}(dt) \\
= \int_{B} \varphi(t, r) \p(\zeta = r|\eta = t)\p_{\eta}(dt).
\end{gather*}
Hence \eqref{V} is valid for $\varphi(t,v)=\ind\{t\in A, v\in E\}$ where one can take
arbitrary $A\in \mathcal{B}(\mathbb{R})$ and $E\subset S$.
Taking into account that any measurable function $\varphi:\mathbb{R}\times S\to \mathbb{R}$
can be approximated by finite linear combinations of the considered functions of the type
$\ind\{A\}\ind\{E\}$ we come to desired statement \eqref{V}. We also note that
$\e |\e (\xi|\eta, \zeta)| \leq \e|\xi| < \infty$. $\square$

The applications of obtained results to the feature selection involving the mutual information
estimation
will be provided in the forthcoming paper.

\vskip0.2cm
{\bf Acknowledgements} The work is supported by the Russian Science Foundation under grant 14-21-00162 and performed at the Steklov Mathematical Institute of Russian Academy of Sciences.
The authors are grateful to Professor V.I.Bogachev for the helpful discussions.

\end{document}